\documentclass[english,10pt]{IEEEtran}
\usepackage{amsmath}
\usepackage{mathrsfs}
\usepackage{amsthm}
\usepackage{amssymb}
\usepackage{graphicx}
\usepackage{cite}
\usepackage{subfigure}
\usepackage{color}

\makeatletter
  \theoremstyle{definition}
  \newtheorem{defn}{\protect\definitionname}
	\theoremstyle{definition}
  
  \theoremstyle{plain}
  \newtheorem{thm}{\protect\theoremname}
  \theoremstyle{plain}
  \newtheorem{lem}{\protect\lemmaname}
  \theoremstyle{plain}
  
  \theoremstyle{remark}

\usepackage{url}
\usepackage{bbm}
\usepackage{epsfig}

\makeatother

\usepackage{babel}
\providecommand{\definitionname}{Definition}
\providecommand{\lemmaname}{Lemma}
\providecommand{\propositionname}{Proposition}
\providecommand{\remarkname}{Remark}
\providecommand{\theoremname}{Theorem}
\providecommand{\examplename}{Example}
\begin{document}
% Change ieeeproof to proof
%If you are using the amsthm package (or one of the ams documentclasses), then the symbol used is given by the command \qedsymbol. It can be redefined:
\renewcommand{\qedsymbol}{$\blacksquare$}

\newcommand{\tcb}[1]{\textcolor{black}{#1}}
\newcommand{\tcba}[1]{\textcolor{black}{#1}}

\title{A Geometric Approach to Aggregate Flexibility Modeling of Thermostatically
Controlled Loads}

\author{Lin Zhao, Wei Zhang, He Hao, and Karan Kalsi
\thanks{
This work was supported in part by the National Science Foundation under
grant ECCS-1309569 and grant CNS-1552838. (\textit{Corresponding author: He Hao.})
\newline
\indent L. Zhao and W. Zhang are with the Dept. of Electrical and Computer Engineering, The Ohio State University, Columbus, OH, USA, 43210
{\small (email: zhao.833@osu.edu; zhang.491@osu.edu)}
\newline 
\indent H. Hao and K. Kalsi are with Pacific Northwest National Laboratory, Richland, WA, USA, 99352
{\small (email: he.hao@pnnl.gov; karanjit.kalsi@pnnl.gov)}}
}
\maketitle
\begin{abstract}
Coordinated aggregation of a large population of thermostatically controlled loads (TCLs) presents a great potential to provide various ancillary services to the grid. One of the key challenges of integrating TCLs into system level operation and control is developing a simple and portable model to accurately capture their aggregate flexibility. In this paper, we propose a geometric approach to model the aggregate flexibility of TCLs. We show that the set of admissible power profiles of an individual TCL is a polytope, and their aggregate flexibility is the Minkowski sum of the individual polytopes. In order to represent their aggregate flexibility in an intuitive way and achieve a tractable approximation, we develop optimization-based algorithms to approximate the polytopes by the homothets of a given convex set. As a special application, this set is chosen as a \emph{virtual battery model} and the corresponding optimal approximations are solved efficiently by equivalent linear programming problems. Numerical results show that our algorithms yield significant improvement in characterizing the aggregate flexibility over existing modeling methods. We also conduct case studies to demonstrate the efficacy of our approaches by coordinating TCLs to track a frequency regulation signal from the Pennsylvania-New Jersey-Maryland (PJM) Interconnection.

%A large population of thermostatically controlled loads (TCLs) can
%be coordinated to provide various ancillary services to the grid.
%Effective operation and coordination of TCLs requires an accurate
%and simple model to capture their aggregate power flexibility. This
%paper proposes a novel geometric approach to model the aggregate flexibility
%of a given TCL population. We adopt a discrete-time formulation of
%the individual TCL dynamics over a finite time-horizon. The set of
%admissible power profiles of an individual TCL is shown to be a polytope.
%The aggregate flexibility is then modeled by the Minkowski sum of
%the individual polytopes. In order to have a tractable approximation
%of the aggregate flexibility, we propose to optimally inner and outer
%approximate the individual polytopes by the homothets of a given convex
%set, repectively. As a special application, this set can be chosen
%as a virtual battery model and the corresponding optimal approximation
%problems can be solved efficiently by equivalent linear programs.
%Our simulations show significant improvement in characterizing the
%aggregate flexibility using virtual battery models over the existing
%modeling methods.
\end{abstract}

\section*{Nomenclature}
\addcontentsline{toc}{section}{Nomenclature}

\begin{IEEEdescription}[\IEEEusemathlabelsep\IEEEsetlabelwidth{$V_1,V_2,V_3$}]

%\item[Thermostatically Controlled Loads (TCLs) Notations:]
\item[$\theta$] Temperature of a TCL system.
\item[$\theta_r$] User specified temperature set-point.
\item[$\delta$] Sampling time of a TCL system.
\item[$\Delta$] Half of the temperature deadband.
\item[$C_{th}$] Thermal capacitance.
\item[$R_{th}$] Thermal resistance.
\item[$P_0$] Nominal power (baseline power) that maintains the set-point temperature.
\item[$P_m$] Rated power of a TCL system.
\item[$q(k)$] Operating state "`OFF/ON"' at time instant $k$.
\item[$\eta$] Coefficient of performance of the power consumption.
\item[$\Omega^i$] The set of the $i^\text{th}$ TCL's parameters.
\item[$\Omega_o$] The set of mean TCL parameters.
%\item[$N$}] The number of the TCLs.}
\item[$\mathcal{T}$] Time horizon set $\{1,2,\dots m\}$.

\item[$A^i$, $B^i$] $m$-dimensional system/input matrix associated with the vector representation of the $i^\text{th}$ TCL's discrete dynamics.
\item[$C^i$] The vector representation of the $i^\text{th}$ TCL's initial condition.
\item[$\Lambda^i$] The inverse of the matrix $A^i$.

%\item[Virtual Battery Notations:]
\item[$a$] Energy dissipation rate.
\item[$X(k)$] Energy state of a virtual battery at time $k$.
\item[$U(k)$] Power supply/draw of a virtual battery at time~$k$.
\item[$C$] Initial condition vector of a virtual battery, $C=[aX(0),0,\cdots,0]^T$.
\item[$\underline{D},\overline{D}$] Lower/Upper power limits of a virtual battery. %, $\underline{D}=[D_-(k)]$, $\overline{D}=[D_+(k)]$.}
\item[$\underline{E},\overline{E}$] Lower/Upper energy capacity limits of a virtual battery. %, $\underline{E}=[E_-(k)]$, $\overline{E}=[E_+(k)]$.}
\item[$\phi$] The set of the virtual battery parameters, $\phi=\{C,\underline{D},\overline{D},\underline{E},\overline{E}\}$.

\item[$\mathbb{B}_{o}$] Prototype virtual battery model.
\item[$\mathbb{B}_{s}$, $\mathbb{B}_{n}$] Sufficient/Necessary virtual battery model.
\item[$\mathbb{B}_{s}^*$, $\mathbb{B}_{n}^*$] Optimal sufficient/necessary virtual battery model.
\item[$\mathbb{B}_{s}^\dagger$, $\mathbb{B}_{n}^\dagger$] Suboptimal sufficient/necessary virtual battery model.
\item[$\mathbb{B}_{s}^\diamond$, $\mathbb{B}_{n}^\diamond$] Sufficient/Necessary virtual battery model obtained using the algorithms proposed in \cite{Hao2015, Hao2015a}.

\item[$\mathcal{P}$] Exact aggregate flexibility.
\item[$\mathcal{P}^{i}$] Individual flexibility of the $i^\text{th}$ TCL.
\item[$\mathcal{P}_{o}$] Prototype set.

\item[$\beta,\beta^i,\beta_*^i$] Scaling factors.
\item[$t,t^i,t_*^i$] Translation vectors.

\item[$I_m$] Identity diagonal matrix of dimension $m$.
\item[$\mbox{diag}(s;-1)$] Lower subdiagonal matrix consisting of element $s$ with appropriate dimension.
\item[$\biguplus$] Operator for Minkowski sum.

\end{IEEEdescription}

\section{Introduction}
Renewable energy resources such as wind and solar have a high degree of variability. Recent studies show that deep penetration of variable generations into the power grid requires substantial reserves from the generation side and flexible consumption via demand response \cite{makarov2009operational, helman2010resource, CallawayPIEEE}. Thermostatically controlled loads (TCLs) such as air conditioners, heat pumps, water heaters, and refrigerators are an important class of demand response assets due to their resource size and inherent flexibility. It is well recognized that coordination of TCLs presents a huge potential to provide various services to the grid, such as frequency regulation, energy arbitrage, renewable integration, and peak shaving, etc. \cite{lu2004state, Callaway2009, CallawayPIEEE, Mathieu2013, Hao2015, Mathieu2015, li2015}. However, to integrate a large population of TCLs into system level operation and control, a fundamental challenge is to construct a simple and user-friendly model to manage them. This model should be able to accurately capture their aggregate flexibility, while being amenable to system level optimization and control.

The existing literature on aggregate modeling of TCLs can be generally divided into two categories: modeling the population dynamics of the loads, and characterizing the set of admissible aggregate power profiles. In the first category, the studies focus on establishing dynamical equations that describe the probability density evolution of a population of TCLs. These include partial differential equations \cite{malhame1985electric, Callaway2009, BashashFathy2013, Zhao2015} and Markov chains \cite{lu2004state, Mathieu2013, Zhang2013, sanandaji2014fast}. However, in order to reproduce the population dynamics accurately, these models often require fine gridding of the state space, which is computationally expensive \cite{Zhang2013, Mathieu2013}. Moreover, the above methods do not explicitly characterize the ex-ante flexibility that TCLs can offer to the grid.

To address the above issue, the second category of aggregate modeling aims to characterize the set of admissible aggregate power profiles that can be consumed by the TCL population without violating any comfort or operational constraint \cite{Hao2015, Hao2015a, Sanandaji2014,  Mathieu2015, BarotTaylor2017}. The set of admissible power profiles represents the \textit{aggregate flexibility} of the TCL population. Such models are for ex-ante planning, which can be easily incorporated into solving various problems, such as multi-period optimal power flow problems~\cite{VrettosAndersson2013,ZhangShenMathieu2016}, or unit commitment problems~\cite{TrovatoTengStrbac2016}, among others. The aggregate flexibility model characterizes the power capacity of the load population, and thus can assist their provision of ancillary services under the system level coordination~\cite{PJM_M12}.

In the literature, the aggregate flexibility is often modeled as a virtual battery model \cite{Hao2015, Hao2015a, Mathieu2015, ZhangShenMathieu2016,HughesDominguezGarciaPoolla2016}. The virtual battery model is a scalar linear system that resembles a simplified battery dynamics parameterized by charge and discharge power limits, energy capacity limits, and self-discharge rate. However, the existing virtual battery models for characterizing the aggregate flexibility are very conservative \cite{Hao2015, Hao2015a}, especially for a TCL population with heterogeneous model parameters. The authors in \cite{Sanandaji2014} proposed several ways to improve the flexibility characterization, but only certain specific battery parameters were optimized independently under special cases of limited TCL population heterogeneities. 

This paper proposes a novel geometric approach which is able to characterize the aggregate flexibility of heterogeneous TCLs more accurately. We show that the power flexibility of an individual TCL can be represented by a polytope, and the \emph{aggregate flexibility} is the Minkowski sum of these polytopes. However, an exact computation of this Minkowski sum is numerically intractable when the number of TCLs is large. Therefore, we estimate the aggregate flexibility by a subset and a superset of the Minkowski sum of the individual flexibilities. Specifically, we first approximate each individual flexibility polytope by its subset and superset respectively, and then calculate the Minkowski sum of the resulting approximations accordingly. The key to facilitating the second step is to restrict the approximation sets to be the homothets (i.e., the dilation and translation) of a given convex set, the latter of which will be referred to as the \emph{prototype set}. Hence, for each TCL, we compute the maximum inner approximation and the minimum outer approximation of its flexibility polytope with respect to the homothets of the prototype set. Moreover, we show if the prototype set is chosen as a polytope, then the optimization problems can be formulated as linear programming problems, and therefore can be solved very efficiently. %We show that the Minkowski sum of the inner/outer approximations is still a homothet of the prototype set and is a sufficient/necessary approximation of the aggregate flexibility.

The above proposed geometric approach provides a general framework for aggregating a large number of constrained linear dynamical systems when the summation of individual quantities is of interest. In particular, the virtual battery modeling in \cite{Hao2015, Hao2015a, Mathieu2015} can be viewed as a special case by choosing the prototype set as the virtual battery model. %We show that the obtained virtual battery model is optimal in terms of set inclusion. Additionally, we develop a suboptimal method for approximating the individual flexibility by taking advantage of the special structure of the virtual battery model. Compared to the optimal method, the suboptimal method substantially reduces the numerical complexity but generally yields more conservative estimation of the aggregate flexibility. 
Compared to the optimization methods proposed in~\cite{Hao2015, Hao2015a, Mathieu2015}, our approach takes advantage of the geometric information of the flexibility polytopes and optimizes over additional decision variables which represent the translation vector. These features improve the modeling accuracy significantly and can deal with much stronger parameter heterogeneity. We show that with $10\%$, $20\%$, and $30\%$ TCL parameter heterogeneities, our approach can improve the flexibility characterization accuracy by as much as $129\%$, $141\%$, and $156\%$ respectively. Moreover, we demonstrate the efficacy of our geometric approach through an example of providing frequency regulation service to the grid, where we control the aggregate power of a population of TCLs to track a regulation signal from the PJM Interconnection \cite{PJMinterconnection}. We show that the proposed approach substantially increases the regulation capacity that TCLs can provide to the ancillary service market, and the dispatched regulation signal can be followed successfully without violating any comfort or operational constraint of TCLs.

Other closely-related works on aggregate flexibility modeling of power system demand-side resources include~\cite{TrangbaekBendtsen2012,MuellerSundstroemSzaboEtAl2015,ZhaoHaoZhang2016, BarotTaylor2017}. These methods are based on the general idea of computing the exact or approximate Minkowski sum of different types of polytopes. An outer approximation of the Minkowski sum of general polytopes was proposed in~\cite{BarotTaylor2017}. However, the number of inequality constraints resulted from this method is non-deterministic and is very large in general. Besides, the outer approximation cannot guarantee the feasibility of the aggregate power profile. The references~\cite{TrangbaekBendtsen2012,MuellerSundstroemSzaboEtAl2015,ZhaoHaoZhang2016}  deals with the so-called resource polytopes, which arise from the flexibility modeling of deferrable loads, such as plug-in electric vehicles (PEVs), dishwashers, among many others. In particular, the exact Minkowski sum of resource polytopes was considered in~\cite{TrangbaekBendtsen2012}. However, it cannot be applied to the flexibility polytopes resulted from TCL systems. In addition, it cannot deal with high dimensional polytopes since the number of inequality constraints increases exponentially with the system dimension. Moreover, an inner approximation of the resource polytopes using Zonotopes was proposed in~\cite{MuellerSundstroemSzaboEtAl2015}, whereas using general polytopes was investigated in~\cite{ZhaoHaoZhang2016} via a lift and projection method.

The rest of the paper is organized as follows. In Section~\ref{sec:modeling}, we present the problem statement. The geometric approach to aggregate flexibility characterization is proposed in Section~\ref{sec:characterization}. In Section~\ref{sec:VB_FC}, we apply the geometric approach to obtain the virtual battery models. We demonstrate its efficacy through numerical examples and case studies in Section~\ref{sec:simulations}. Finally, we summarize our research and discuss the future work in Section~\ref{sec:conclusions}.

\section{Modeling of TCL and Flexibility}\label{sec:modeling}
In this section, we first present a nonlinear switching model that governs the temperature dynamics of a TCL. To facilitate aggregate flexibility modeling, we adopt a constrained linear system model to approximate the power consumption of the switching model. Based on this linear system model, we define the individual and aggregate flexibility of TCLs. It is worth mentioning that the linear system model is only employed for analysis purpose, and the nonlinear switching model is used in all the simulation studies presented in Section~\ref{sec:simulations}. 

%we first introduce a switching system TCL model, which will be used in the simulation of providing frequency regulation. To facilitate our analysis of aggregate flexibility modeling, we then consider a constrained linear system model. Finally, the concept of the aggregate flexibility is defined based on the second model.

%Before we proceed, we present the notations used in this paper. 
%\textbf{Notation}: A polytope $\mathcal{P}$ is a solution set of a system of finite linear inequalities: $\mathcal{P}:=\{x:Ax\leq c\}$, where  $\leq$ (or $<$, $\geq$,$>$) denotes elementwise operation. A polytope $\mathcal{P}\subset\mathbb{R}^{m}$ is called full dimensional if it contains an interior point in $\mathbb{R}^{m}$. We use $\biguplus$ to denote the Minkowski sum of multiple sets, and $\oplus$ of two sets. Additionally, ${\bf 1}_{m}$ represents the $m$-dimensional row vector of all ones, and $I_{m}$ represents $m$-dimensional identity matrix. We use $[u(k)]$ to denote the vector whose $k^\text{th}$ element is $u(k)$. %\textbf{The index set $\{1,2,\cdots,N\}$ is denoted by $\mathcal{N}$. For two matrices $u$ and $v$ with the same column dimension, we write $(u,v)$ for $[u^{T},v^{T}]^{T}$ when no confusion arises.}
\subsection{Nonlinear Switching Model of TCLs}

The temperature evolution of a TCL can be described by a discrete-time switching model \cite{Callaway2009,Mathieu2013,Hao2015}:
\begin{equation}\label{eq:TCLCtsDyn}
\theta(k)=
a\theta(k-1)+(1-a)(\theta_{a}-b q(k) P_{m}),
\end{equation}
where $\theta(k)$ is the TCL temperature at time step $k$, $\theta_{a}$ is the ambient temperature whose dynamics are much slower than $\theta(k)$, $P_{m}$ is the rated power, and $q(k)\in \{0,1\}$ is a binary variable representing the operating state \textquotedblleft OFF/ON\textquotedblright{} of the system. The model parameters $a$ and $b$ are related to the thermal capacitance $C_{th}$, thermal resistance $R_{th}$, and coefficient of performance $\eta$ of the system by $a=e^{-\Delta T/(R_{th}C_{th})}\approx 1- \Delta T/(R_{th}C_{th})$ and $b=R_{th} \eta$, where $\Delta T$ is the sampling time. Without loss of generality, we assume each TCL is a cooling device with $P_{m}>0$. The TCL switches between \textquotedblleft ON\textquotedblright{} and \textquotedblleft OFF\textquotedblright{} subject to the following local control rules,
\begin{equation}
q(k)=\begin{cases}
1, & \theta(k-1)\geq\theta_{r}+\Delta,\\
0, & \theta(k-1)\leq\theta_{r}-\Delta,\\
q(k-1), & \mbox{otherwise},
\end{cases}\label{eq:localcontrol}
\end{equation}
where $\theta_{r}$ is the user-specified temperature set-point and $\Delta>0$ is half of the deadband.
\vspace*{-0.2cm}
\subsection{Linear System Model of TCLs}
To aggregate the flexibility of TCLs, the above switching model \eqref{eq:TCLCtsDyn}-\eqref{eq:localcontrol} is very challenging for analysis due to its nonlinearity. Therefore, we consider a linear system model to approximate it,  
\begin{equation}\label{eq:discretize}
\theta(k)=a\theta(k-1)+(1-a)(\theta_{a}-b P(k)),
\end{equation}
where $P(k)\in[0,P_m]$ is a \emph{continuous} variable instead of a binary input of $\{0,P_m\}$. It is shown in \cite{Hao2015, Hao2015a, Sanandaji2016} that the aggregate behavior of a large population of TCLs with model \eqref{eq:TCLCtsDyn}-\eqref{eq:localcontrol} can be accurately approximated by model \eqref{eq:discretize}. The continuous power input $P(k)$ can be considered as the average of the binary power input of model \eqref{eq:TCLCtsDyn}-\eqref{eq:localcontrol} over time. Additionally, for a large population of TCLs, the aggregate power of the linear system models can match that of the nonlinear switching models closely. After a change of variables, $x(k)=C_{th}(\theta_{r}-\theta(k))/\eta$, and $u(k)=P(k)-P_{0}(k)$, where $P_{0}(k)=(\theta_{a}-\theta_{r})/b$ is the nominal power that keeps the temperature of model \eqref{eq:discretize} at its set-point, we can rewrite the above model as,
\begin{equation}\label{eq:battery}
x(k)=a x(k-1)+ u(k) \delta,
\end{equation}
where $\delta=(1-a)R_{th}C_{th}\approx \Delta T$. Additionally, the model has energy constraint $x(k)\in[-x_{-},x_{+}]$ with $x_{+}=x_{-}=C_{th}\Delta/\eta$, and input constraint $u(k)\in[-u_{-}(k),u_{+}(k)]$ with $u_{-}(k)=P_{0}(k)$ and $u_{+}(k)=P_{m}-P_{0}(k)$.
%Similar linear system models have also been employed in~\cite{BagagioloBauso2014},~\cite{KizilkaleMalhame2013} and~\cite{LiZhangLianEtAl2016} to facilitate their analysis. Our analysis for aggregate flexibility modeling will be based on (\ref{eq:discretize}).

\vspace*{-0.1cm}
\subsection{Modeling of Flexibility}
We consider a heterogeneous population of $N$ TCLs modeled by~(\ref{eq:battery}). Each TCL is parameterized by $\Omega^{i}:=\{R_{th}^{i},C_{th}^{i},\theta_{r}^{i},\Delta^{i},\eta^{i},\theta^{i}(0),P_{m}^{i}\}$. The aggregate power consumption of a population of TCLs has many feasible solutions that respect all the temperature and power constraints of TCLs. The key of nondisruptive control of TCLs for demand response is to accurately characterize their aggregate power flexibility over a considered time horizon $\mathcal{T}:=\{1,2,\dots m\}$. Before we proceed, we first define the individual and the aggregate flexibilities of TCLs. 

\begin{defn}\label{def:flexibility}
For each TCL $i=1,\cdots,N$, its \emph{individual flexibility} is defined as the set of all admissible power profiles
{\small
\begin{align*}\label{eq:Pi}
\mathcal{P}^{i}=\left\{[u^{i}(k)]\in\mathbb{R}^{m} \Bigg|
\begin{matrix}
x^i(k)=a^i x^i(k-1)+u^i(k)\delta, \ \forall \ k\in\mathcal{T}\\
-u^i_{-}(k)\leq u^i(k)\leq u^i_{+}(k),\ \forall \ k\in\mathcal{T}\\
-x^i_{-}\leq x^i(k)\leq x_{+}^i,\ \forall \ k\in\mathcal{T}
\end{matrix}
\right\}.
\end{align*}
}

\noindent  where $[u^{i}(k)]$ denotes a vector whose $k^\text{th}$ element is $u^{i}(k)$. The \textit{aggregate flexibility} of a population of TCLs is a set of power profiles satisfying
\[
\mathcal{P}=\left\{ U\in\mathbb{R}^{m} \Big| U=\sum_{i=1}^{N}u^{i},\forall \ u^{i}\in\mathcal{P}^{i}\right\} .
\]
The aggregate flexibility can be written as 
\[
\mathcal{P}=\biguplus_{i=1}^{N}\mathcal{P}^{i},
\] 
where $\biguplus$ denotes the Minkowski sum. 
\end{defn}

The set $\mathcal{P}$ contains all the aggregate power profiles that are admissible to the population of TCLs. However, the expression of set $\mathcal{P}$ is very abstract, and it is challenging to integrate it into the power system level operation and control. To represent the aggregate flexibility in an intuitive way, we define a simple and portable virtual battery model which will be used to describe the aggregate flexibility of TCLs.

\begin{defn}\label{def:virtualbattery}
An $m$-horizon discrete-time virtual battery model is a set of power
profiles satisfying
{\small
\begin{align*}
\left\{\left[U(k)\right] \in \mathbb{R}^m \Bigg|
\begin{matrix}
X(k)=a X(k-1)+ U(k)\delta,\ \forall \ k\in\mathcal{T}\\
-D_{-}(k)\leq U(k)\leq D_{+}(k),\ \forall \ k\in\mathcal{T}\\
-E_{-}(k)\leq X(k)\leq E_{+}(k),\ \forall \ k\in\mathcal{T}
\end{matrix}
\right\}.
\end{align*}
}
The virtual battery model is specified by parameters 
\[
\phi:=\left(a,X(0),D_{-}(k),D_{+}(k),E_{-}(k),E_{+}(k),\forall \ k\in\mathcal{T}\right),
\]
and we write it compactly as $\mathbb{B}(\phi)$. In addition, $\mathbb{B}(\phi)$ is called sufficient if $\mathbb{B}(\phi)\subset\mathcal{P}$, and it is called necessary if $\mathbb{B}(\phi)\supset\mathcal{P}$. 
\end{defn}
Note that the definition of the virtual battery model comes naturally from the definition of the individual flexibility. We can regard $U(k)$ as the power draw of the battery and $X(k)$ as its charging state which indicates the level of the energy stored in the battery. The quantities $D_{+}(k)$ and $D_{-}(k)$ represent the time-varying charging/discharging rate limits, and $E_{+}$ and $E_{-}$ represent upper/lower energy capacity limits relative to the nominal energy level of the battery. In addition, the parameter $a$ represents the self-discharge rate of the battery, which is due to the thermal exchange between the inner air and the ambient environment. We will show in the next section that the virtual battery model offers us great convenience in describing and characterizing the aggregate flexibility of TCLs.

\section{Geometric Approach to Flexibility Characterization}\label{sec:characterization}
In this section, we present a geometric interpretation of the aggregate flexibility. Additionally, we show that it is generally intractable to compute the exact aggregate flexibility. Therefore, \tcba{optimal approximations of the aggregate flexibility are proposed. They are further formulated as linear programming problems which can be solved very efficiently.}%we formulate its maximum inner approximation (sufficient characterization) and minimum outer approximation (necessary characterization) problems. They are further converted to equivalent linear programming problems.}%In particular, these sets can be chosen to have a battery interpretation. We then develop two methods to optimally approximate the aggregate flexibility using virtual battery models.

\subsection{Polytope Interpretation of Flexibility}
A polytope $\mathcal{Q}$ is a solution set of a system of finite linear inequalities: $\mathcal{Q}:=\{U \in \mathbb{R}^m|FU\leq H\}$, where $\leq$ denotes elementwise inequality. A polytope $\mathcal{Q}\subset\mathbb{R}^{m}$ is called full dimensional if it contains an interior point in $\mathbb{R}^{m}$. In this subsection, we show that the aggregate flexibility  of TCLs can be represented by a polytope. Denoting the state and input vectors by $X^i=[x^i(k)]$ and $U^i=[u^i(k)]$, we rewrite (\ref{eq:battery}) as, 
\begin{equation}
A^iX^i=B^iU^i+C^i,\label{eq:sysMa}
\end{equation}
where $A^i=I_{m}+\mbox{diag}(-a^i;-1)$ is a lower bidiagonal matrix with $1$'s on the main diagonal, and $-a^i$'s on the lower subdiagonal, $B^i=\delta I_{m}$, in which $I_{m}$ denotes the $m$-dimensional identity matrix, and $C^i=[a^i x^i(0),0,\cdots,0]^{T}$. \tcba{The inverse of $A^i$ can be derived in an explicit form with polynomials of $a^i$ as its entries.} It will be denoted by $\varLambda^i$ in the sequel. Additionally, the constraint sets for $X^i$ and $U^i$ are 
\begin{equation}
-\underline{U}^i\leq U^i\leq\bar{U}^i,\ -\underline{X}^i\leq X^i\leq\bar{X}^i,\label{eq:constraints}
\end{equation}
where $\underline{U}^i=[u^i_{-}(k)],$ $\bar{U}^i=[u^i_{+}(k)],$ $\underline{X}^i={\bf 1}_{m}x^i_{-},$ and $\bar{X}^i={\bf 1}_{m}x^i_{+}$, in which ${\bf 1}_{m}$ is the $m$-dimensional column vector of all ones. 

Using (\ref{eq:sysMa}) and (\ref{eq:constraints}), the individual flexibility $\mathcal{P}^{i}$ of the $i^\text{th}$ TCL can be expressed as 
\begin{equation}
\mathcal{P}^{i}=\mathcal{U}^{i}\cap\mathcal{X}^{i},\label{eq:DisAdmissible}
\end{equation}
where
\begin{gather*}
\mathcal{U}^{i}=\left\{ U^{i}\in\mathbb{R}^{m}|-\underline{U}^{i}\leq U^{i}\leq\bar{U}^{i}\right\},\\
\mathcal{X}^{i}=\left\{ U^{i}\in\mathbb{R}^{m}|-\underline{X}^{i}\leq\varLambda^{i}B^{i}U^{i}+\varLambda^{i}C^{i}\leq\bar{X}^{i}\right\}.
\end{gather*}
Since $-\infty<-u_{-}^{i}<u_{+}^{i}<+\infty$, it is straightforward to show that $\mathcal{U}^{i}$ is a full dimensional hyper rectangular (and thus a polytope). Similarly, because $\varLambda^{i}$ and $B^{i}$ are invertible and $-\infty<-x_{-}^{i}<x_{+}^{i}<+\infty$, we can show that $\mathcal{X}^{i}$ is also a full dimensional polytope. It then follows from \cite{Henk1997} that their intersection $\mathcal{P}^{i}=\mathcal{U}^{i}\cap\mathcal{X}^{i}$ is a polytope if $\mathcal{P}^{i}\neq\emptyset$. Moreover, it can be proven that the Minkowski sum of polytopes $\mathcal{P}=\biguplus_{i=1}^{N}\mathcal{P}^{i}$ is also a polytope \cite{Schneider1993}.

%\begin{thm}
%\label{thm:polytope}$\forall i\in\mathcal{N}$, the individual flexibility
%$\mathcal{P}^{i}$ of discrete-time finite time horizon admissible
%power profiles is a polytope, and thus the aggregate flexibility $\mathcal{P}=\biguplus_{i=1}^{N}\mathcal{P}^{i}$
%is also a polytope.
%\end{thm}

For each TCL, its individual flexibility $\mathcal{P}^{i}$ can be determined by \eqref{eq:sysMa}-\eqref{eq:DisAdmissible}. However, the numerical complexity of calculating their Minkowski sum is prohibitively expensive when the number of TCLs is large. In fact, calculating the Minkowski sum of two sets $\mathcal{Q}_{1}$ and $\mathcal{Q}_{2}$ when they are polytopes specified by facets is NP-hard since the facets of the obtained polytope can grow exponentially with the number of the facets of $\mathcal{Q}_{1}$ and $\mathcal{Q}_{2}$~\cite{Tiwary2008, Weibel2007}. Therefore, we take an alternative route and find its maximum inner (subset) approximation and minimum outer (superset)  approximation instead.

\subsection{Optimal Approximations of the Aggregate Flexibility}
%Instead of directly calculating the exact set $\mathcal{P}$, we aim to find the optimal inner and outer approximations of it. 
In this subsection, we aim to find sets $\mathcal{P}_{s}$ and $\mathcal{P}_{n}$ such that $\mathcal{P}_{s}\subset\mathcal{P}\subset\mathcal{P}_{n}$. Any such sets $\mathcal{P}_{s}$ and $\mathcal{P}_{n}$ will be referred to as the \emph{sufficient} approximation and the \emph{necessary} approximation, respectively. Given a power profile $U$, if $U\notin\mathcal{P}_{n}$, we can conclude that $U$ is not an admissible aggregate power profile for TCLs. On the other hand, if $U\in\mathcal{P}_{s}$, then there exists a decomposition of $U$  such that $U = \sum_{i=1}^NU^{i}$, and $U^{i}\in\mathcal{P}^{i}$ for all $i =1, \cdots, N$.

\begin{figure}[tb]
\begin{center}
\includegraphics[width=0.9\linewidth]{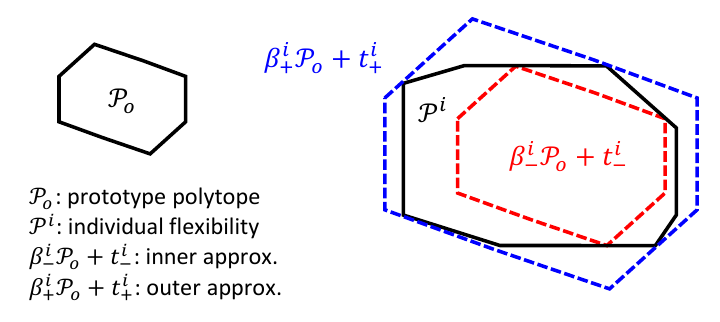}
\vspace{-5pt}
\caption{Optimal inner and outer approximations of the individual flexibility polytope $\mathcal{P}^{i}$ with respect to the prototype set $\mathcal{P}_{o}$}\label{fig:InnerOuter}
\end{center}
\vspace{-10pt}
\end{figure}

Given a compact convex set $\mathcal{P}_{o}$, we call $\beta^{i}\mathcal{P}_{o}+t^{i}:=\{U^i|U^i=\beta^{i}\xi+t^{i},\forall\ \xi\in\mathcal{P}_{o}\}$ a \textit{homothet} of $\mathcal{P}_{o}$, \tcba{that is, the dilation and translation of $\mathcal{P}_{o}$,} where $\beta^{i}>0$ is a scaling factor and $t^{i}\in\mathbb{R}^{m}$ is a translation factor. Since all $\mathcal{P}^{i}$'s have the same structure \eqref{eq:DisAdmissible}, \tcba{we conduct the inner and outer approximations of each $\mathcal{P}^{i}$ with respect to a given set, $\mathcal{P}_{o}$. The set $\mathcal{P}_{o}$ will be referred to as the \emph{prototype set} hereafter. Specifically, we will find within the homothets of $\mathcal{P}_{o}$ the optimal approximations of each $\mathcal{P}^{i}$'s.} Fig.~\ref{fig:InnerOuter} illustrates this idea using a 2-dimensional example (i.e., the time horizon is taken as $m=2$), where $\beta^{i}_{+}\mathcal{P}_{o}+t^{i}_{+}$ is a minimum outer approximation of $\mathcal{P}^{i}$ and $\beta^{i}_{-}\mathcal{P}_{o}+t^{i}_{-}$ is a maximum inner approximation of $\mathcal{P}^{i}$.

The benefit of using homothets to approximate $\mathcal{P}^{i}$'s is that it admits an efficient calculation of their Minkowski sum. It was shown in~\cite{Schneider1993} that given a convex set $\mathcal{Q}$, non-negative scalars $\beta^{1}$ and $\beta^{2}$, and any scalars $t^{1}$ and $t^{2}$, the calculation of their Minkowski sum is simply as follows,
\begin{equation}
\tcba{(\beta^{1}\mathcal{Q}+t^1)\biguplus(\beta^{2}\mathcal{Q}+t^2)=(\beta^{1}+\beta^{2})\mathcal{Q} +(t^1+t^2), \label{eq:HomothetMink}}
\end{equation}
which says that the Minkowski sum of homothets of a convex set reduces to the sum of the scaling factors and the sum of the translation factors. \tcba{Therefore, once we have the $\mathcal{P}_{o}$-homothetic approximations of the individual flexibility polytopes, the approximation of the aggregate flexibility can be calculated very easily.} 

\tcba{Our focus now becomes how to choose the $\mathcal{P}_{o}$-homothet that optimally approximates $\mathcal{P}^{i}$. Specifically, for each $i=1,\cdots, N$, we are interested in finding the maximal $\mathcal{P}_{o}$-homothet that is contained in $\mathcal{P}^{i}$, which can be mathematically expressed as,}
\begin{equation}
\begin{array}{ll}
\underset{\beta^{i},t^{i}}{\mbox{maximize}} & \beta^{i}\\
\mbox{subject to:} & \beta^{i}\mathcal{P}_{o}+t^{i}\subset\mathcal{P}^{i},\\
 & \beta^{i}>0,
\end{array}\label{eq:Inner}
\end{equation}
\tcba{and the minimal $\mathcal{P}_{o}$-homothet that contains $\mathcal{P}^{i}$,}
\begin{equation}
\begin{array}{ll}
\underset{\beta^{i},t^{i}}{\mbox{minimize}} & \beta^{i}\\
\mbox{subject to:} & \beta^{i}\mathcal{P}_{o}+t^{i}\supset\mathcal{P}^{i},\\
 & \beta^{i}>0.
\end{array}\label{eq:Outer}
\end{equation}
\tcba{In the above optimization problems}, the optimality is in the sense of inclusion, i.e., if $(\beta^{i}_{*},t^{i}_{*})$ is an optimal solution, then there is no other homothet of $\mathcal{P}_{o}$ contained in between $\beta^{i}_{*}\mathcal{P}_{o}+t^{i}_{*}$ and $\mathcal{P}^{i}$.
\tcba{We will refer to the optimal solution of \eqref{eq:Inner} and \eqref{eq:Outer} as the Maximum Inner Approximation (MIA) and Minimum Outer Approximation (MOA), respectively.} For convenience, we denote the solutions of problems \eqref{eq:Inner} and \eqref{eq:Outer} by $(\beta^{i},t^{i})=\mbox{MIA}(\mathcal{P}^{i},\mathcal{P}_{o})$ and $(\beta^{i},t^{i})=\mbox{MOA}(\mathcal{P}^{i},\mathcal{P}_{o})$, respectively.

\tcba{In order to have tractable solutions of the MIA and MOA problems, we have to specify $\mathcal{P}_{o}$ beforehand.} Since each $\mathcal{P}^{i}$ is a polytope, \tcba{to achieve a good approximation of it,} it is natural to choose $\mathcal{P}_{o}$ as a polytope too. Furthermore, we show that the MIA and MOA problems under such choice can be solved efficiently by equivalent linear programming problems \cite{Eaves1982}. The specific optimization algorithms to solve for $(\beta_{i},t_{i})=\mbox{MIA}(\mathcal{P}^{i},\mathcal{P}_{o})$ and $(\beta_{i},t_{i})=\mbox{MOA}(\mathcal{P}^{i},\mathcal{P}_{o})$ are derived as follows. In view of~(\ref{eq:sysMa})-(\ref{eq:DisAdmissible}), the individual
flexibility polytope can be written as $\mathcal{P}^{i}=\left\{ U^i\in\mathbb{R}^{m} |F^{i}U^i\leq H^{i}\right\}$, where
\begin{align*}
F^{i} & =(I_{m},-I_{m},\Lambda^{i}B^{i},-\Lambda^{i}B^{i}), \\
H^{i} & =(\bar{U}^{i},\underline{U}^{i},\bar{X}^{i}-\Lambda^{i}C^{i},\underline{X}^{i}+\Lambda^{i}C^{i}),
\end{align*}
where $(x,y)$ denotes the matrix $[x^T,y^T]^T$ for two matrices $x$ and $y$ with the same number of columns.
Moreover, if $\mathcal{P}_{o}$ is chosen to have the form $\{U\in\mathbb{R}^{m}:FU\leq H\}$, then we have the following theorem:
\begin{thm}\label{thm:AlgInner}
The optimal solution of $\mbox{MIA}(\mathcal{P}^{i},\mathcal{P}_{o})$ is given by $\beta^{i}_{*}=1/s^{i}_*$, and $t^{i}_{*}=-r^{i}_*/s^{i}_*$, if $(s^{i}_*, r^{i}_*, G_*)$ is an optimal solution of the following linear programming problem,
\begin{equation}
\begin{array}{ll}
\underset{s^{i}>0,G\geq0,r^{i}}{\mathrm{minimize}} & s^{i}\\
\mathrm{subject\ to:} & GF=F^{i},\\
 & GH\leq s^{i}H^{i}+F^{i}r^{i}.
\end{array}\label{eq:LPOP1}
\end{equation}
Similarly, the optimal solution $(\beta^i_*,t^i_*)=\mbox{MOA}(\mathcal{P}^{i},\mathcal{P}_{o})$ is solved by
\begin{equation}
\begin{array}{ll}
\underset{\beta^{i}>0,G\geq0,t^{i}}{\mathrm{minimize}} & \beta^{i}\\
\mathrm{subject\ to:} & GF^{i}=F,\\
 & GH^i\leq \beta^{i}H+Ft^{i}.
\end{array}\label{eq:LPOP2}
\end{equation}
\end{thm}
\begin{IEEEproof}
See Appendix \ref{subsec:Proof_Prop}.
\end{IEEEproof}

\tcba{The above theorem provides algorithms for solving the optimal inner/outer $\mathcal{P}_{o}$-homothetic approximations of individual flexibility polytopes. Furthermore, it is easy to see that if $\mathcal{P}^i\subset\mathcal{Q}^i$ for some polytope $\mathcal{Q}^i$ for each $i$, then we also have $\biguplus_i\mathcal{P}^i\subset\biguplus_i\mathcal{Q}^i$. Therefore, the inner/outer approximation of the aggregate flexibility can be obtained by the Minkowski sum of the obtained $\mathcal{P}_{o}$-homothets. By the formula \eqref{eq:HomothetMink}, we see that $\biguplus_i(\beta^i\mathcal{P}_{o}+t^i)=(\sum_i\beta^i)\mathcal{P}_{o}+\sum_i{t^i}$ is a sufficient approximation of $\mathcal{P}$ if $(\beta^{i},t^{i})=\mbox{MIA}(\mathcal{P}^{i},\mathcal{P}_{o})$, and is a necessary approximation of $\mathcal{P}$ if $(\beta^{i},t^{i})=\mbox{MOA}(\mathcal{P}^{i},\mathcal{P}_{o})$. Note that the obtained aggregate flexibility approximations are also $\mathcal{P}_{o}$-homothets. Therefore, some desired properties of the aggregate flexibility model can be achieved through specifying $\mathcal{P}_{o}$'s structure. In the next section, we will choose $\mathcal{P}_{o}$ as virtual battery models and derive the battery parameters for the approximated aggregate flexibility. We further propose computationally more efficient virtual battery modeling methods by exploiting their special structures.}

%As we discussed, the set $\mathcal{P}_{o}$ is better chosen to have the particular structure induced from the constrained model~(\ref{eq:discretize}). Such a structure can be interpreted naturally by a simple battery model. Therefore, in the next subsection, we will define this special class of polytope as the virtual battery models. Its structure will also help us to further develop different methods to approximate $\mathcal{P}^{i}$.

\section{Virtual Battery based Flexibility Characterization}\label{sec:VB_FC}
We now consider a special case of our geometric approach proposed in the last section, where the prototype set $\mathcal{P}_o$ is chosen as the virtual battery model (see Definition \ref{def:virtualbattery}). In the sequel, we assume that the battery parameter $a$ is predetermined and fixed (e.g., taking the mean of all TCL parameters $a^{i}$'s), and focus on estimating its power limits, energy capacity, and initial energy state. Hence, we denote $\phi=(C,\underline{D},\bar{D},\underline{E},\bar{E})$ as the parameter of the virtual battery in the sequel, where the vector notations are given by $C=[a X(0),0,\cdots,0]^{T}$, $\underline{D}=\left[D_{-}(k)\right]$, $\bar{D}=\left[D_{+}(k)\right]$, $\underline{E}=\left[E_{-}(k)\right]$, and $\bar{E}=\left[E_{+}(k)\right]$. We will use $\mathbb{B}_{*}$ as the short notation for the virtual battery $\mathbb{B}(\phi_{*})$, where $*$ denotes the subscript $o$, $s$, and $n$ meaning prototype, sufficient, and necessary virtual battery, respectively.

\begin{figure}[tb]
\begin{center}
\includegraphics[clip,width=.85\linewidth]{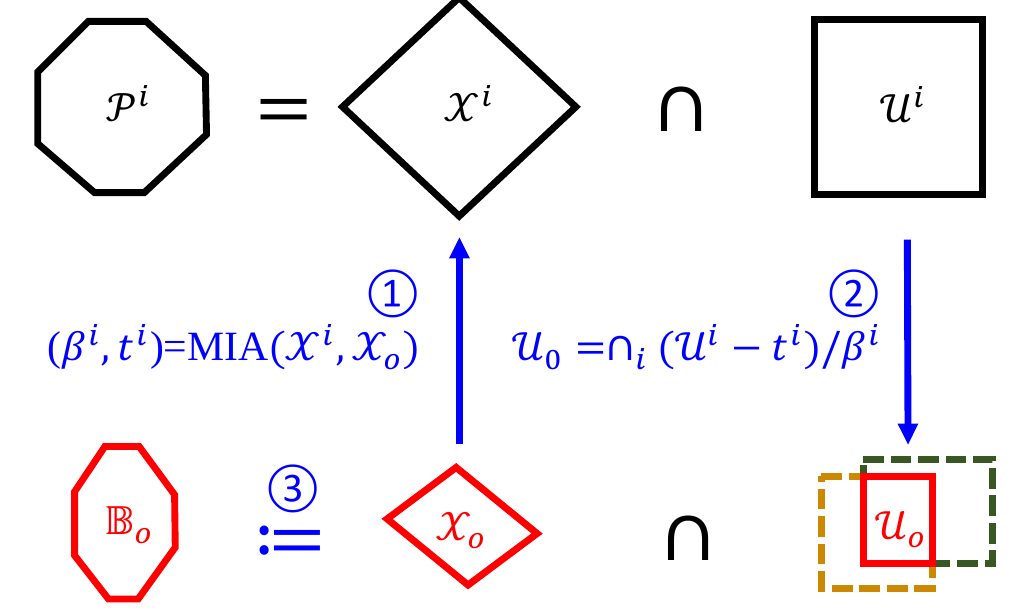}
\vspace{-5pt}
\caption{\tcba{Diagram of the suboptimal sufficient virtual battery characterization}}\label{fig:suboptimal_Scheme}
\end{center}
\vspace{-8pt}
\end{figure}

\subsection{Optimal Virtual Battery}
In this subsection, we choose $\mathcal{P}_{o}$ to be in the same form as $\mathcal{P}^i$'s with parameter $\Omega_{o}$ being the mean of all TCL parameters $\Omega^{i}$'s. \tcba{Clearly, the $\mathcal{P}_{o}$ of this choice is also a virtual battery model. To emphasize it,  we will denote this prototype virtual battery by $\mathbb{B}_{o}$,} with parameters $\phi_{o}=(C_{o},\underline{D}_{o},\bar{D}_{o},\underline{E}_{o},\bar{E}_{o})$ obtained from the average of the TCL parameters. In view of~(\ref{eq:sysMa})-(\ref{eq:DisAdmissible}), the prototype virtual battery can be expressed as $\mathbb{B}_{o}=\left\{ U\in\mathbb{R}^{m}:\ FU\leq H\right\} $, where $F=(I_{m},-I_{m},\Lambda B,-\Lambda B)$ and $H=(\bar{D}_{o},\underline{D}_{o},\bar{E}_{o}-\Lambda C_{o},\underline{E}_{o}+\Lambda C_{o})$. As discussed in the previous section, we can approximate each $\mathcal{P}^{i}$ by $\beta^{i}\mathbb{B}_{o}+t^{i}$, where $\beta^i$ and $t^i$ are the optimal solutions of \eqref{eq:LPOP1} or \eqref{eq:LPOP2}. Now let $\beta=\sum_{i=1}^{N}\beta^{i}$ and $t=\sum_{i=1}^{N}t^{i}$.
We will show that $\beta\mathbb{B}_{o}+t$ is a sufficient or necessary
battery depending on whether $\{\beta^{i},t^{i}\}$'s are obtained from \eqref{eq:LPOP1} or \eqref{eq:LPOP2}. 
\begin{thm}
\label{thm:SuffVB}$\mathbb{B}_{*}:=\beta\mathbb{B}_{o}+t$ is a sufficient
battery if $(\beta^{i},t^{i})=\mbox{MIA}(\mathcal{P}^{i},\mathbb{B}_{o}),$
and is a necessary battery if $(\beta^{i},t^{i})=\mbox{MOA}(\mathcal{P}^{i},\mathbb{B}_{o})$.
The battery parameter $\phi_{*}=(C,\underline{D},\bar{D},\underline{E},\bar{E})$
is given by
\begin{equation}
\begin{cases}
C=\beta C_{o},\\
\underline{D}=\beta\underline{D}_{o}-t, & \bar{D}=\beta\bar{D}_{o}+t,\\
\underline{E}=\beta\underline{E}_{o}-\Lambda Bt, & \bar{E}=\beta\bar{E}_{o}+\Lambda Bt.
\end{cases}\label{eq:Suff1bounds}
\end{equation}
In addition, $\forall \ U\in\mathbb{B}_{s}$, the individual admissible
power profile is given by
\begin{equation}
U^{i}=\frac{\beta^{i}}{\beta}(U-t)+t^{i},\ \forall \ i.\label{eq:Ui}
\end{equation}
\end{thm}
\begin{IEEEproof}
See Appendix~\ref{subsec:Proof_Theorem2}.
\end{IEEEproof}

%The set $\mathbb{B}_{s}$ (or $\mathbb{B}_{n}$) will be referred to as the optimal sufficient (or optimal necessary) virtual battery model with respective to $\mathbb{B}_{o}$.
\tcba{This theorem is a direct result of the proposed general aggregate flexibility modeling method by choosing $\mathcal{P}_o$ as the virtual battery model. The advantage of the virtual battery model lies in that the resulted aggregate flexibility model has a very simple form, which is very desirable in practice for system level optimization. In addition, as will be shown later through simulations, it can give a very accurate approximation of the exact aggregate flexibility.}
Moreover, we emphasize that $\beta^{i}$'s and $t^{i}$'s for different $i$'s can be computed in parallel. This makes our algorithm scalable and it can be executed efficiently even with a large number of TCLs.
\subsection{Suboptimal Virtual Battery}

To further reduce the computational complexity, we will propose a suboptimal method. In the above optimal method, we approximate $\mathcal{P}^{i}=\mathcal{U}^{i}\cap\mathcal{X}^{i}$ with respect to $\mathbb{B}_{o}$ as a whole. Alternatively, we can approximate its components $\mathcal{U}^{i}$ and $\mathcal{X}^{i}$ separately. Although this is generally more conservative than the optimal method, it substantially reduces the numerical complexities since the corresponding linear programming problem has a smaller dimension.

We first consider the sufficient battery characterization. For notational convenience, we denote $\mathbb{B}_{o}=\mathcal{U}_{o}\cap\mathcal{X}_{o}$ with the same structure of~(\ref{eq:DisAdmissible}). There are mainly 3 steps of this method, which are illustrated in Fig.~\ref{fig:suboptimal_Scheme} by a 2-dimensional example (i.e., the time horizons $m=2$). Assuming $\mathcal{X}_{o}$ is given as one of the components of $\mathbb{B}_{o}$, in the first step each $\mathcal{X}^{i}$ is approximated with respect to $\mathcal{X}_{o}$, which yields $(\beta^{i},t^{i})=\mbox{MIA}({\mathcal{X}}^{i},\mathcal{X}_{o})$. \tcba{Next, to guarantee the homothetic transformation between $\mathbb{B}_{o}$ and $\mathcal{P}^{i}$, the $\mathcal{U}_{o}$-component of $\mathbb{B}_{o}$ can be determined by enforcing $\beta^i\mathcal{U}_{o}+t^i\subset\mathcal{U}^{i}$ for all $i$, or equivalently by enforcing}
\begin{equation}
\mathcal{U}_{o}\subset\cap_i(\mathcal{U}^{i}-t^i)/\beta^i. \label{eq:cal_U0} 
\end{equation}
Since $\mathcal{U}^{i}$'s are all hyper rectangulars, the right hand side of \eqref{eq:cal_U0} can be calculated exactly which is given below in~\eqref{eq:U0powerbds}. Clearly, this is the largest $\mathcal{U}_{o}$ that satisfies \eqref{eq:cal_U0}. The last step is to obtain $\mathbb{B}_{o}$ as $\mathcal{X}_{o}\cap\mathcal{U}_{o}$, and then $\mathbb{B}_{s}$ is obtained by $\beta\mathbb{B}_{o}+t$.

We next prove that with $\mathbb{B}_{o}$ obtained in this way, $\beta\mathbb{B}_{o}+t$ is a sufficient battery.

\begin{thm}
\label{thm:SuffVB_V2} Let $(\beta^{i},t^{i})=\mbox{MIA}(\mathcal{X}^{i},\mathcal{X}_{o})$
and $\mathcal{U}_{o}$ be given by $\mathcal{U}_{o}=\{U\in\mathbb{R}^{m}|-\underline{U}\leq U\leq\bar{U}\}$, where
\begin{equation}
\underline{U}=\underset{i}{\min}\frac{\underline{U}^{i}+t^{i}}{\beta^{i}},\ \bar{U}=\underset{i}{\min}\frac{\bar{U}^{i}-t^{i}}{\beta^{i}},\label{eq:U0powerbds}
\end{equation}
and $\min$ is element-wise. Then $\mathbb{B}_{s}:=\beta\mathbb{B}_{o}+t$ is a sufficient battery. The battery parameters
$\phi_{s}$ are given by~(\ref{eq:Suff1bounds}). Additionally, $\forall \ U\in\mathbb{B}_{s}$,
the individual admissible power profile can be obtained by~(\ref{eq:Ui}).
\end{thm}
\begin{IEEEproof}
See Appendix \ref{subsec:Proof_Theorem3}.
\end{IEEEproof}

It is worth mentioning that Theorem \ref{thm:SuffVB_V2} can be easily adapted to obtain a suboptimal necessary battery \cite{Zhao2016}. However, its estimation of the power limits can be very inaccurate. In the following theorem, we develop a suboptimal necessary battery modeling method that generalizes the methods proposed in~\cite{Hao2015,Hao2015a} to improve the power limits characterization. The main idea is that the two components of the necessary battery can be obtained independently. \tcba{In particular, its power limits can be obtained by simply adding up the individual power limits, while the energy capacity limits are obtained through the $\mathcal{X}_o$-homothetic approximation of each $\mathcal{X}^i$.}

\begin{thm}
\label{thm:NeceVB_V2} Let $(\beta_{i},t_{i})=\mbox{MOA}(\mathcal{X}^{i},\mathcal{X}_{o})$
and $\mathcal{U}_{o}=\uplus_{i}\mathcal{U}^{i}$. Then $\mathbb{B}_{n}:=(\beta\mathcal{X}_{o}+t)\cap\mathcal{U}_{o}$
is a necessary battery. The battery parameters are given by
\begin{equation}
\begin{cases}
C=\beta C_{o},\\
\underline{D}=\sum_{i=1}^{N}\underline{U}^{i}, & \bar{D}=\sum_{i=1}^{N}\bar{U}^{i},\\
\underline{E}=\beta\underline{E}_{o}-\Lambda Bt, & \bar{E}=\beta\bar{E}_{o}+\Lambda Bt.
\end{cases}\label{eq:NeceBounds2nd}
\end{equation}
\end{thm}

\begin{IEEEproof}
See Appendix \ref{subsec:Proof_SuboptNece}.
\end{IEEEproof}

\tcba{Our approach (Theorems~\ref{thm:SuffVB}-\ref{thm:NeceVB_V2}) includes the approach proposed in~\cite{Hao2015, Hao2015a} as a special case. As an example, it can be seen from~(\ref{eq:NeceBounds2nd}) that the suboptimal necessary battery calculates the power bounds in the same way as~\cite{Hao2015,Hao2015a}, i.e., by summing over the individual power bounds. In addition, it optimizes the approximation of $\mathcal{X}^i$ using $\mathcal{X}_o$-homothet. In fact, if we further drop these optimization schemes (i.e., MIA($\mathcal{X}^{i}$,$\mathcal{X}_{o}$) and MOA($\mathcal{X}^{i}$,$\mathcal{X}_{o}$) in Theorems~\ref{thm:SuffVB_V2}-\ref{thm:NeceVB_V2}, and instead obtain the energy bounds from~(\ref{eq:sysMa}) using matrix norm inequalities, then we will get the counterpart of
the modeling method proposed in~\cite{Hao2015,Hao2015a} in the discrete-time
finite-horizon case. Moreover, the optimal method proposed in Theorem~\ref{thm:SuffVB} approximates $\mathcal{P}^{i}$ directly, and therefore offers the best performance among all these methods.}

The major computation of the suboptimal methods only involves solving $\mbox{MIA}(\mathcal{X}^{i}, \mathcal{X}_{o})$ or $\mbox{MOA}(\mathcal{X}^{i}, \mathcal{X}_{o})$, which involves about $50\%$ fewer constraints and
$75\%$ fewer decision variables as compared to the optimal method. Therefore, it can be solved much faster than the optimal case.

\begin{figure}
\centering
\subfigure[Power Limits]{\includegraphics[clip,width=\columnwidth]{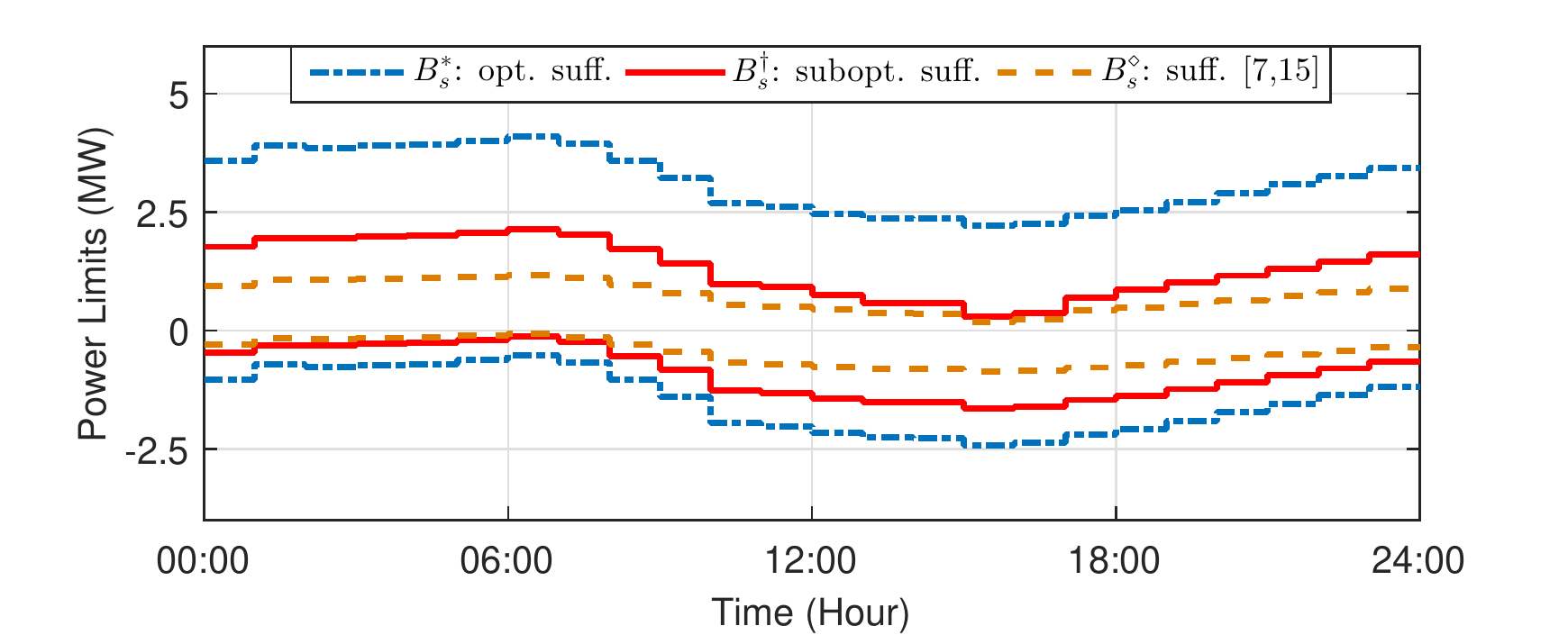}}
\subfigure[Energy Limits]{\includegraphics[clip,width=\columnwidth]{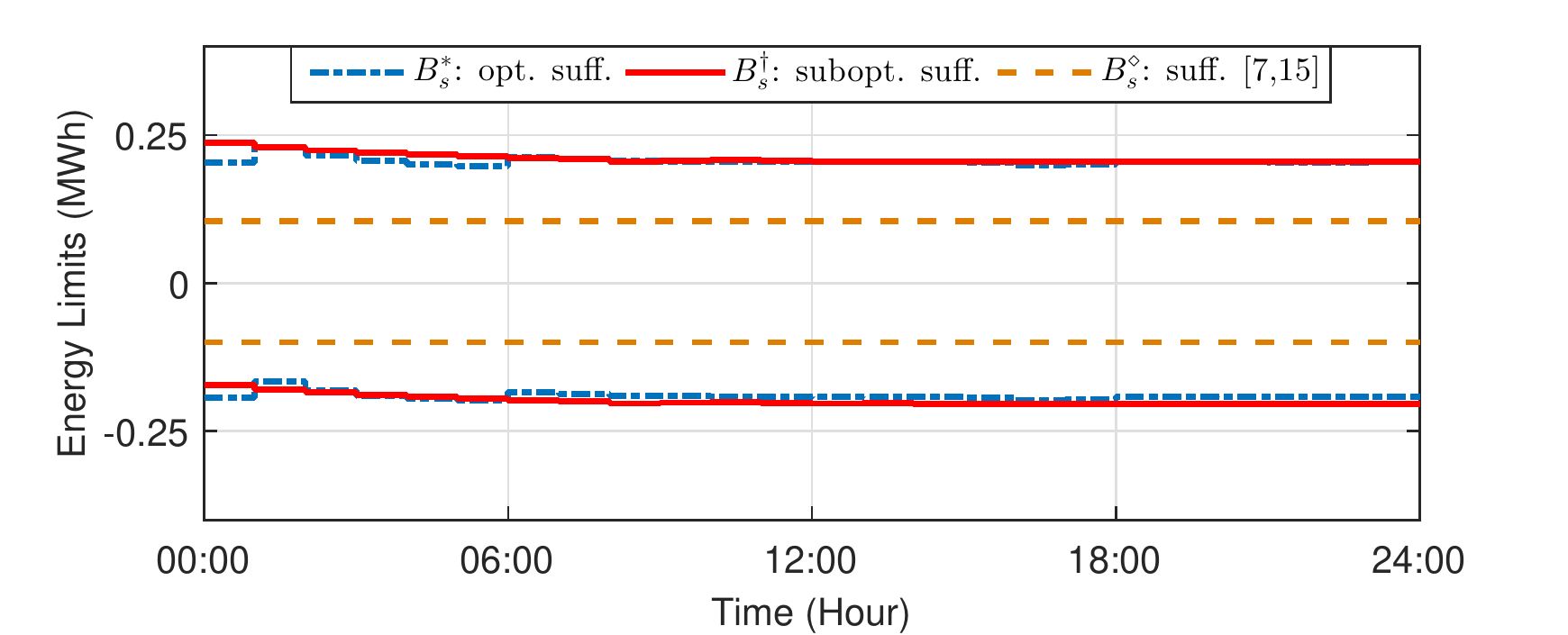}}
\caption{Sufficient battery comparison among $\mathbb{B}_{s}^{*}$, $\mathbb{B}_{s}^\dagger$, and $\mathbb{B}_{s}^{\diamond}$.}
\label{fig:SufficientBattery}
\end{figure}

\section{Case Studies}\label{sec:simulations}
In this section, we first compare the characterized flexibilities using our geometric approach and the method in \cite{Hao2015,Hao2015a}. We next demonstrate the efficacy of our approach through an example of providing frequency regulation service. We then show that the population of TCLs using our flexibility characterization approach can provide more regulation capacity to the grid while achieving excellent tracking of the regulation signal and guaranteeing the thermal comfort of the end users. 

We consider a population of 1000 heterogeneous TCLs. Their model parameters  $\Omega^i$'s are assumed to be uniformly distributed \cite{Callaway2009, Mathieu2013, Hao2015}, e.g., the thermal capacities $C_{th}$'s are uniformly distributed within $[(1-\epsilon)\bar{C}_{th},(1+\epsilon)\bar{C}_{th}]$, where $\bar{C}_{th}$ is the mean value, and $\epsilon$ models the degree of heterogeneity. Additionally, the ambient temperature profile is picked as a typical hot summer day in Columbus, OH \cite{Weather}. 

\tcba{For the comparison of the virtual batteries,} note that if a sufficient battery has both larger power and energy limits than another one, we claim that the former extracts more flexibility than the latter and the latter is more conservative. On the other hand, if a necessary battery has both larger power and energy limits than another one, we claim that the former overestimates more flexibility than the latter, and the latter is more accurate.

\subsection{Performance Comparison}

We first calculate the optimal sufficient battery using Theorem~\ref{thm:SuffVB}. The MIA and MOA problems are solved using the GLPK linear programming solver~\cite{glpk2006} interfaced with YALMIP~\cite{Lofberg2004}.  The blue dash-dot lines in Fig. \ref{fig:SufficientBattery} (a) (respectively, (b)) are the lower and upper power (respectively, energy) limits of the optimal sufficient battery (denoted by $\mathbb{B}_{s}^{*}$). If a given power profile $U\in \mathbb{R}^m$ belongs to $\mathbb{B}_{s}^{*}$, then it must lie between the two blue dash-dot lines in Fig.~\ref{fig:SufficientBattery} (a), and the energy state resulting from this power profile through the battery dynamical equation $X(k)=aX(k-1)+\delta U(k)$ (see Definition \ref{def:virtualbattery}) must lie between the blue dash-dot lines in Fig. \ref{fig:SufficientBattery} (b). In other words, given a power profile $U$ which lies between the two blue dash-dot lines in Figs. \ref{fig:SufficientBattery} (a), if the associated energy state vector $X$ also lies between the blue dash-dot lines in Figs. \ref{fig:SufficientBattery} (b), then $U$ is a feasible power trajectory of the TCL population.

Our first observation is that the ambient temperature has a significant impact on the aggregate flexibility. For example, when the temperature is the lowest around 6:00 AM (the temperature profile is not shown here due the space limit), the baseline aggregate power consumption is the lowest at this time. As a result, the lower power limit (which corresponds to the largest possible down regulation capacity) is the smallest (see Fig. 3), while the upper power limit (which corresponds to the largest possible up regulation capacity) is the largest. The situation is inverted when the temperature reaches the highest around 3 PM, since the nominal aggregate power consumption is the highest.

Next, we use Theorem \ref{thm:SuffVB_V2} to obtain the suboptimal sufficient battery (denoted by $\mathbb{B}_{s}^\dagger$) and the red solid lines in Fig.~\ref{fig:SufficientBattery} (a) (respectively, (b)) are respectively its lower and upper power (respectively, energy) limits of the suboptimal sufficient battery. Roughly speaking, the optimal sufficient battery $\mathbb{B}_{s}^{*}$ extracts more flexibility than the suboptimal battery $\mathbb{B}_{s}^\dagger$, since $\mathbb{B}_{s}^{*}$ has larger power limits than $\mathbb{B}_{s}^\dagger$ and their energy limits are similar. Moreover, we compare our geometric approach with the characterization methods in \cite{Hao2015,Hao2015a}, where the sufficient battery (denoted by $\mathbb{B}_{s}^\diamond$) is obtained by solving the optimization problem in \cite[Theorem 3]{Hao2015a}. The orange dashed lines in Figs.~\ref{fig:SufficientBattery} (a) and (b) represent $\mathbb{B}_{s}^\diamond$. It can be seen that both $\mathbb{B}_{s}^{*}$ and $\mathbb{B}_{s}^\dagger$ extracts more flexibility than $\mathbb{B}_{s}^\diamond$, i.e., $ \mathbb{B}_{s}^{*}\supset\mathbb{B}_{s}^\diamond $ and $\mathbb{B}_{s}^\dagger \supset \mathbb{B}_{s}^\diamond$, since both their power and energy limits are larger than those of $\mathbb{B}_{s}^\diamond$.

%For the sufficient virtual battery, both methods can be viewed as calculating Minkowski sum of \textquoteleft similar\textquoteright{} objects. The latter method also has the freedom to optimize over the scaling factors $\{\beta^{i},\forall i\in\mathcal{N}\}$. Theorem 3 of \cite{Hao2015a} proposed a linear programming optimization problem to find the Pareto-set of the power and energy bounds $\rho:=(\underline{D},\bar{D},\underline{E},\bar{E})$ of the battery. In the case of symmetric bounds $\bar{D}=\underline{D}$ or $\underline{E}=\bar{E}$, the optimal bounds can be calculated analytically (see~\cite[Table III]{Hao2015}). 

\begin{figure}[tb]
\centering
\subfigure[Power Limits]{\includegraphics[clip,width=\columnwidth]{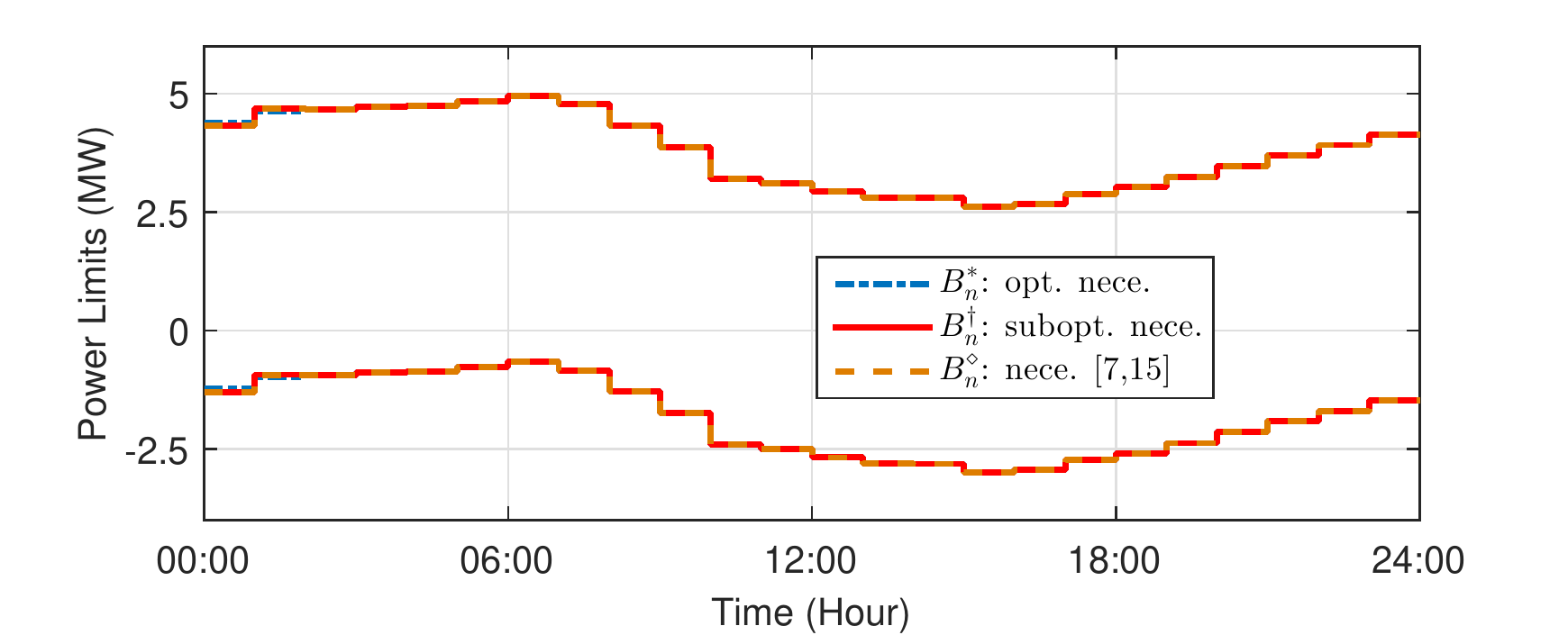}}
\subfigure[Energy Limits]{\includegraphics[clip,width=\columnwidth]{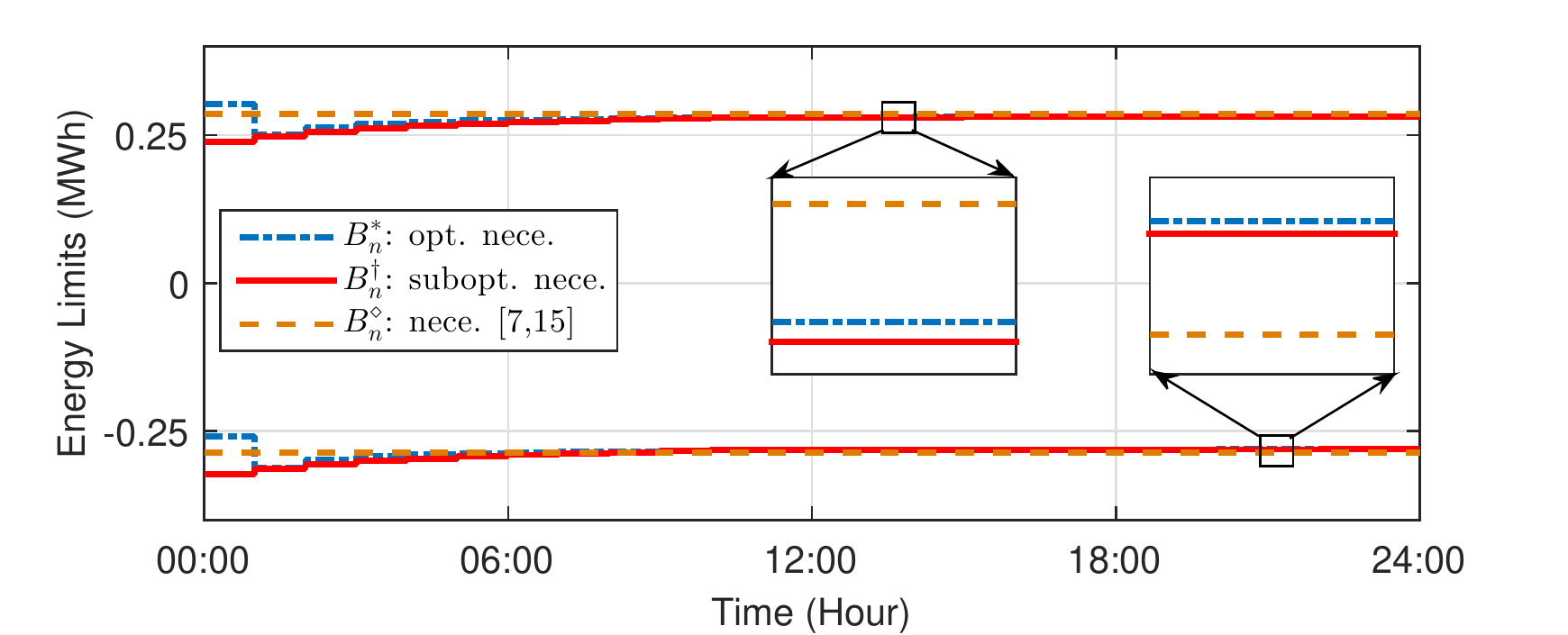}}
\caption{Necessary battery comparison among $\mathbb{B}_{n}^{*}$,
$\mathbb{B}_{n}^{\dagger}$, and $\mathbb{B}_{n}^\diamond.$}
\label{fig:NecessaryBattery}
\end{figure}

%In addition, we plot the necessary battery proposed in~\cite{Hao2015} (denoted by $\hat{\mathbb{B}}_{n}$) in Figs.~\ref{fig:Power242_3Methods}-\ref{fig:Energy242_3Methods} by the blue dotted lines. It is interesting to observe that by comparing $\mathbb{B}_{s}^{*}$ with $\hat{\mathbb{B}}_{n}$, the exact aggregate flexibility can be approximated even more accurately, since it must satisfy $\mathbb{B}_{s}^{*}\subset\mathcal{P}\subset\hat{\mathbb{B}}_{n}$ while $\mathbb{B}_{s}^{*}$ and $\hat{\mathbb{B}}_{n}$ are very close to each other in terms of the distances between their power and energy bounds.

Furthermore, we compare in Fig. \ref{fig:NecessaryBattery} the power and energy limits of the optimal necessary battery $\mathbb{B}_{n}^{*}$, suboptimal necessary battery $\mathbb{B}_{n}^\dagger$, and necessary battery $\mathbb{B}_{n}^\diamond$ obtained in \cite{Hao2015,Hao2015a}. In obtaining the optimal necessary battery, we choose the suboptimal necessary battery as its prototype battery model. It can be seen all the necessary batteries have similar estimations since their power and energy limits are similar. Even though no strict inclusion relationship is present in our numerical comparison, $\mathbb{B}_{n}^{*}$ and  $\mathbb{B}_{n}^{\dagger}$ are generally more accurate than $\mathbb{B}_{n}^\diamond$, since their energy limits are slightly tighter than $\mathbb{B}_{n}^\diamond$ most of the time, as shown in Fig. \ref{fig:NecessaryBattery} (b).
%In fact, our simulation results show that the obtained optimal necessary battery differs from the suboptimal one only by a translation factor, while having the scaling factor 1 with respect to the suboptimal one. This implies that the optimal necessary battery is not unique, and the same is the MOA problems.

Compared to the method in \cite{Hao2015, Hao2015a},  our approach takes advantage of the geometric information of each individual flexibility set, and thus improves the approximation of the aggregate flexibility. In addition, the flexibility characterization method in \cite{Hao2015, Hao2015a} requires the power and energy bounds (e.g., $u_-(k), u_+(k), x_-$, and $x_+$) to be non-negative for each TCL. This non-negativity requirement restricts the allowable degree of parameter heterogeneity. In contrast, the proposed geometric approach removes such a restriction, and allows us to characterize the aggregate flexibility of a population of TCLs where their model parameters are strongly heterogeneous. We show the performance improvement of the proposed approach at different heterogeneity degrees $\epsilon$ in Table \ref{tab:hetero}. For the convenience of comparison, we assume only the thermal parameters $C_{th}$ and $R_{th}$ are heterogeneous. The numbers in the table represent the average percentage improvement of the power and energy limits as compared to the battery $\mathbb{B}_{s}^\diamond$ and $\mathbb{B}_{n}^\diamond$ in \cite{Hao2015, Hao2015a}, which are calculated as
\begin{align*}
0.5\left|\Gamma(D,D^{\diamond})+\Gamma(E,E^{\diamond})\right|,
\end{align*}
where $\Gamma(D,D^{\diamond}):=\frac{1}{m}\sum_{k=1}^{m}\frac{D_{-}(k)+D_{+}(k)-D_{-}^{\diamond}(k)-D_{+}^{\diamond}(k)}{D_{-}^{\diamond}(k)+D_{+}^{\diamond}(k)}$, 
$D_{-}(k)$, $D_{+}(k)$, $\forall \ k\in\mathcal{T}$ denote the power limits of the optimal or suboptimal batteries, \tcba{and the ones with superscript $\diamond$ represent the power limits of the virtual batteries proposed in \cite{Hao2015, Hao2015a}. The other term $\Gamma(E,E^{\diamond})$ is defined similarly for the energy bounds}. It can be seen from Table \ref{tab:hetero} that the stronger the heterogeneity is, the larger the improvement can be achieved by the proposed approach. %In addition, we see that the suboptimal necessary battery is as good as the optimal necessary battery under the index \eqref{eq:metric}.
\begin{table}
\caption{Improvement of flexibility characterizations under different amounts of heterogeneity.\label{tab:hetero}}
\begin{center}
    \begin{tabular}{| c || c | c || c | c |}
    \hline
   Heterogeneity  & $\mathbb{B}_{s}^{*}$ & $\mathbb{B}_{s}^\dagger$ & $\mathbb{B}_{n}^{*}$ & $\mathbb{B}_{n}^\dagger$ \\ \hline  \hline
   $\epsilon=10\%$ & 129.16\% & 86.13\%   & 0.34\% & 0.34\%  \\ \hline
   $\epsilon=20\%$ & 141.05\% & 87.86\%   & 0.61\% & 0.61\%   \\ \hline
   $\epsilon=30\%$ & 155.94\% & 88.43\%   & 0.82\% & 0.82\%   \\ \hline
    \end{tabular}
\end{center}
\end{table}

\subsection{Providing Frequency Regulation Service}

In order to provide frequency regulation service, each service provider needs to bid its regulation capacities into the day-ahead or hour-ahead ancillary service market. After the market is cleared, each awarded regulating resource will be dispatched a regulation signal $r(t)$ in real-time. The regulation signal will be within the bidded capacity, and it is generally broadcast every 4 seconds depending on the independent system operators. Each regulation resource is obliged to follow this regulation signal accurately since the tracking accuracy will be reflected in the financial settlement. 

\begin{figure}[tb]
\centering
\subfigure[Power trajectory]{\includegraphics[clip,width=\columnwidth]{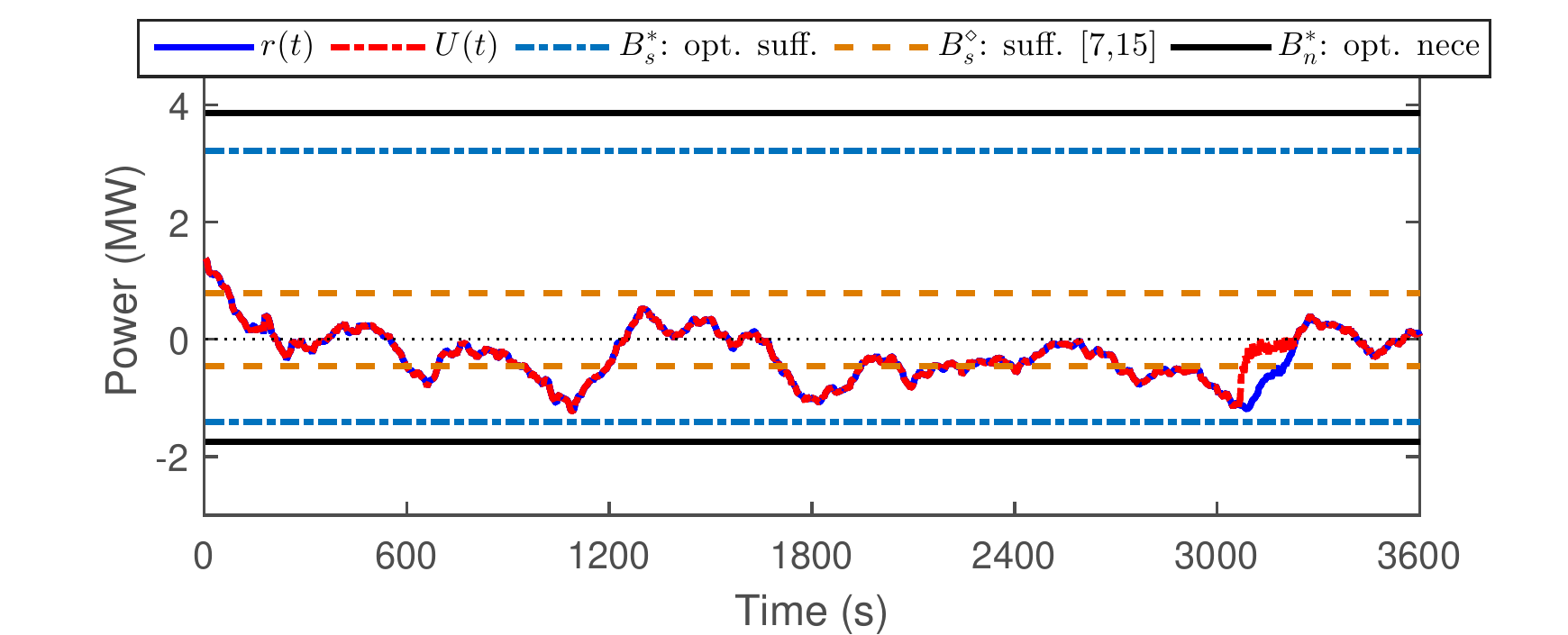}}
\subfigure[Energy state trajectory]{\includegraphics[clip,width=\columnwidth]{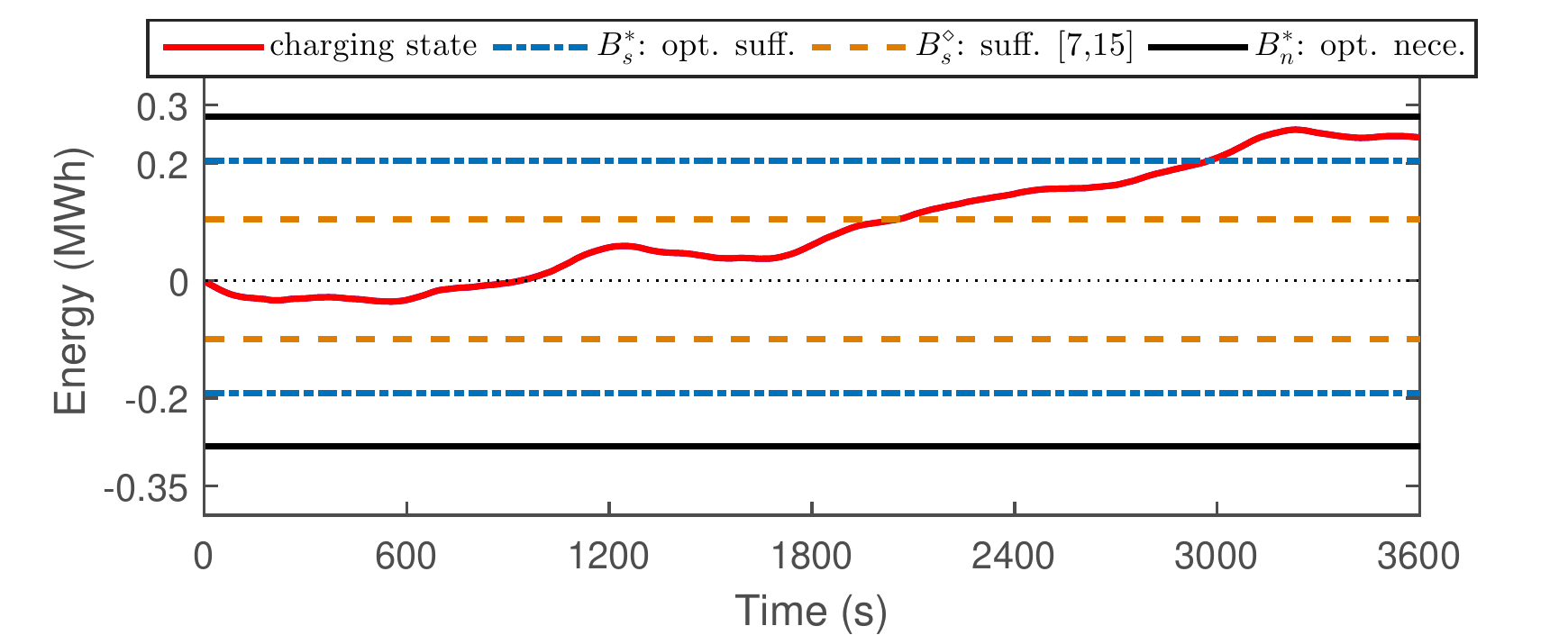}}
\caption{Control of TCLs to track a frequency regulation signal from PJM.}
\label{fig:frequencyregulation}
\end{figure}

Different from the generation resources which has no energy limits, the TCL aggregation has both limits on its power and energy capacity. As a result, the regulation capacity of the TCL aggregation is time-varying, which is coupled over time through the thermal dynamics \eqref{eq:battery}. The virtual battery model gives a characterization of the dynamical changes of the regulation capacity of the TCL aggregation via mimicking the battery dynamics (see Definition \ref{def:virtualbattery}). For example, the power limits of the sufficient virtual battery can be used as a guidance to determine the regulation capacity of the TCL aggregation. Clearly, a less conservative sufficient virtual battery can help improve the estimation of the regulation capacity, and thus increase the potential revenue from providing the regulation service.

To demonstrate coordination of TCLs for frequency regulation,
we first use the proposed geometric approach to
characterize the aggregate flexibility for the population of
TCLs. The lower and upper power limits (dash-dot straight lines in Fig.~\ref{fig:frequencyregulation} (a)) of the proposed optimal sufficient batteries are utilized to determine the regulation capacity of the TCL aggregation. A regulation signal from PJM Interconnection \cite{PJMinterconnection} is then scaled within the power limits of the sufficient batteries. The nonlinear switching model \eqref{eq:TCLCtsDyn}-\eqref{eq:localcontrol} which is sampled at 4-second interval is used in our simulations. We then control the ON/OFF status of the TCLs using the priority-stack-based controller proposed in~\cite{Hao2015} so that the aggregate power of the TCLs minus their baseline tracks the regulation signal. Note that our simulation only involves the lower level control of the TCLs following the
regulation signal, while the estimation of the regulation capacity and the clearing of the
market are assumed to be already done for the real time control. Therefore, we ignored the
network and the associated network constraints at the real time control level.

\begin{figure}[tb]
\begin{center}
\includegraphics[width=\linewidth]{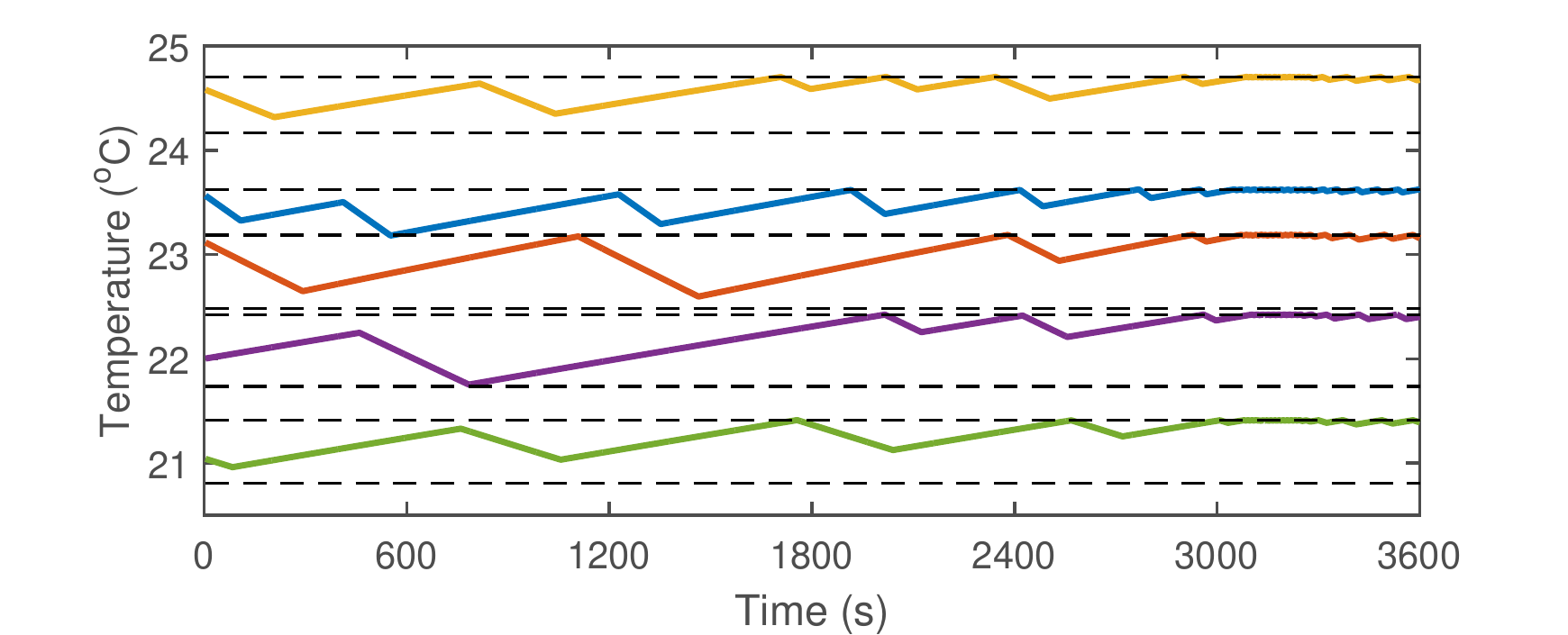}
\caption{Sample TCL temperature trajectories under coordination.}\label{fig:TCLSampleTraj}
\end{center}
\vspace{-10pt}
\end{figure}

The power and energy limits of $\mathbb{B}_{s}^{*}$, $\mathbb{B}_{s}^\diamond$, and $\mathbb{B}_{n}^*$ as well as the regulation signal $r(t)$, the aggregate power deviation from baseline $U(t)$, and the energy state of the virtual battery $X(t)$ are plotted in Fig.~\ref{fig:frequencyregulation}. It can be seen that as long as $r(t)\in\mathbb{B}_{s}^{*}$ (i.e., both the regulation signal $r(t)$ and the resulted charging state are within the power limits and the energy limits of $\mathbb{B}_{s}^{*}$, respectively), $U(t)$ can track $r(t)$ successfully, even when $r(t)$ violates the power and energy limits of the sufficient battery $\mathbb{B}_{s}^\diamond$ obtained in \cite{Hao2015}. This shows that our characterization method is more accurate in estimating the aggregate flexibility of TCLs. Moreover, we observe that $U(t)$ fails to track $r(t)$ shortly after the energy state exceeds the upper energy limit of optimal sufficient battery model $\mathbb{B}_{s}^{*}$ at around 3000 seconds. This again implies that our geometric approach makes a very accurate approximation to the aggregate flexibility of TCLs.

Moreover, Fig.~\ref{fig:TCLSampleTraj} shows the temperature evolutions of several randomly chosen TCLs, where the black dashed lines represent the corresponding allowable temperature bands for each TCLs. We observe that these TCLs are well regulated within the user-specified temperature bands, which means the thermal comfort of end users is strictly respected.

%Therefore, the proposed virtual battery modeling method can be used in various power optimization problems that involve the aggregation of the TCLs to facilitate the real-time computation of near-optimal solutions. 
\section{Conclusions and Future Work}\label{sec:conclusions}

In this paper, we proposed a novel geometric approach to characterize the aggregate flexibility of a population of TCLs. We showed that the power flexibility of an individual TCL could be modeled as a polytope, and their aggregate flexibility was represented by the Minkowski sum of the aforementioned polytopes. However, an exact computation of the Minkowski sum was numerically expensive. We thus developed two optimization-based algorithms to approximate the aggregate flexibility using the maximum inner approximation and minimum outer approximation with respect to the homothets of a prototype set. Additionally, we showed that if the prototype was chosen to be a virtual battery model, our geometric approach extracted more flexibility than existing algorithms in the literature. Moreover, we demonstrated the efficacy of our method through a case study of controlling TCLs to provide regulation service to the grid. We showed that our method could enable TCLs to bid more regulation capacities to the ancillary service market, while achieving excellent tracking of the regulation signal and respecting the thermal comfort requirement of end users. In the future, we are interested in examining the impact of no-short-cycling or minimum off-time constraint on the aggregate flexibility and the aggregate power ramping rate of TCLs.

\appendix{}

\subsection{Proof of Theorem \ref{thm:AlgInner} \label{subsec:Proof_Prop}}

To prove the theorem, we first state the following version of Farkas's lemma \cite{Eaves1982, Mangasarian2002}, which assists in deriving the algorithms for solving the MIA and MOA problems.
\begin{lem}[Farkas' lemma]
\label{lem:Farkas} Suppose that the system of inequalities
$Lx\leq b,$ $L\in\mathbb{R}^{m\times n}$ has a solution and that
every solution satisfies $Mx\leq d,$ $M\in\mathbb{R}^{k\times n}$.
Then there exists $G\in\mathbb{R}^{k\times m}$ , $G\geq0$, such
that $GL=M$ and $Gb\leq d$. The converse is also true.
\end{lem}

After a change of variables $\beta^{i}=1/s^{i},$ and $t^{i}=-r^{i}/s^{i}$, where $s^{i}>0$, the optimization problem (\ref{eq:Inner}) is equivalent to
\begin{equation}
\begin{array}{ll}
\underset{s^{i}>0,r^{i}}{\mbox{minimize}} & s^{i}\\
\text{subject to}: & \mathcal{P}_{o}\subset s^{i}\mathcal{P}^{i}+r^{i}.
\end{array}\label{eq:LPOP1EQV}
\end{equation}
Moreover, we see that $s^{i}\mathcal{P}^{i}+r^{i}$ is the solution
set of $F^{i}U^i\leq s^{i}H^{i}+F^{i}r^{i}$ with respect to $U^i\in\mathbb{R}^{m}$.
Therefore, by Lemma~\ref{lem:Farkas}, there exists a matrix $G$ such that $GF=F^i$, and $GH\leq s^{i}H^{i}+F^{i}r^{i}$. As such, we showed that the optimization problem \eqref{eq:LPOP1EQV} or \eqref{eq:Inner} is equivalent to the linear programming problem \eqref{eq:LPOP1}. Similarly, we can prove that the optimization problem \eqref{eq:Outer} is equivalent to the linear programming problem \eqref{eq:LPOP2}. \hfill \qed

\subsection{Proof of Theorem \ref{thm:SuffVB} \label{subsec:Proof_Theorem2}}

In this section, we only prove the results for the sufficient battery characterization. The necessary
battery characterization can be proved analogously. Suppose $\{\beta^{i},t^{i}\ |\ \forall \ i=1, \cdots, N\}$
are the solutions of the corresponding MIA problems and let $\mathbb{B}_{s}:=\biguplus_{i=1}^{N}(\beta^{i}\mathbb{B}_{o}+t^{i})$.
It can be shown that
\begin{equation}
\mathbb{B}_{s}=\biguplus_{i=1}^{N}\beta^{i}\mathbb{B}_{o}+t=\beta\mathbb{B}_{o}+t,\label{eq:LU0}
\end{equation}
where the last equality is the Minkowski sum of the homothets of $\mathbb{B}_{o}$ \tcba{(recall the formula (\ref{eq:HomothetMink}))}. 
Furthermore, since $\beta^{i}\mathbb{B}_{o}+t^{i}\subset\mathcal{P}^{i}$,
we have
\[
\biguplus_{i=1}^{N}(\beta^{i}\mathbb{B}_{o}+t^{i})\subset\biguplus_{i=1}^{N}\mathcal{P}^{i}=\mathcal{P}.
\]
This implies $\mathbb{B}_{s}$ is a sufficient battery. Its battery parameters $\phi_{s}$ can be obtained by formulating the constraint sets on $U$ from $(U-t)/\beta\in\mathbb{B}_{o}$.
Now $\forall \ U\in\mathbb{B}_{s}$, we obtain from \eqref{eq:LU0} that $(U-t)/\beta\in\mathbb{B}_{o}$. It follows by the scaling and translating of $\mathbb{B}_{o}$
that $\beta^{i}(U-t)/\beta+t^{i}\in\mathcal{P}^{i}$. This completes the proof.
\hfill \qed

\subsection{Proof of Theorem \ref{thm:SuffVB_V2}\label{subsec:Proof_Theorem3}}

From the solution of the $\mbox{MIA}$ problem,
we have $\forall \ i =1, \cdots, N$, $\mathcal{X}^{i}\supset\beta^{i}\mathcal{X}_{o}+t^{i}$.
We next show that $\mathcal{U}^{i}\supset\beta^{i}\mathcal{U}_{o}+t^{i}$,
$\forall \ i =1, \cdots, N$. Since  $\mathcal{U}^{i}$'s
are hyper rectangulars, the largest $\mathcal{U}_{o}$ we can obtain
is $\mathcal{U}_{o}=\cap_{i}\left(\mathcal{U}^{i}-t^{i}\right)/\beta^{i}$. This is equivalent to having the power limits of $\mathbb{B}_{o}$
in~(\ref{eq:U0powerbds}). The above two inclusion relationships
yield
\begin{align*}
\mathcal{P}=\biguplus_{i=1}^{N}\mathcal{X}^{i}\cap\mathcal{U}^{i}&\supset\biguplus_{i=1}^{N}(\beta^{i}\mathcal{X}_{o}+t^{i})\cap(\beta^{i}\mathcal{U}_{o}+t^{i})\\
&=\biguplus_{i=1}^{N}\beta^{i}(\mathcal{X}_{o}\cap\mathcal{U}_{o})+t^{i}=\beta\mathbb{B}_{o}+t,
\end{align*}
where $\mathbb{B}_{o}:=\mathcal{X}_{o}\cap\mathcal{U}_{o}$ and the
last equality is the Minkowski sum of the homothets of $\mathcal{X}_{o}\cap\mathcal{U}_{o}$. It is straightforward to see that $\mathcal{P}\supset\mathbb{B}_{s}:=\beta\mathbb{B}_{o}+t$,
and thus $\mathbb{B}_{s}$ is sufficient. The rest of the proof regarding
the derivation of the parameter $\phi_{s}$ and the power profile
decomposition is the same as those in the proof of Theorem~\ref{thm:SuffVB}.\hfill \qed

\subsection{Proof of Theorem \ref{thm:NeceVB_V2}\label{subsec:Proof_SuboptNece}}
For arbitrary subsets $\mathcal{Q}_{i}$,
$\mathcal{S}_{i}$ of $\mathbb{R}^{m}$,  \tcba{it is easy to verify that} the following
holds,
\[
\biguplus_{i=1}^{N}\left(\mathcal{Q}_{i}\cap\mathcal{S}_{i}\right)\subset\left(\biguplus_{i=1}^{N}\mathcal{Q}_{i}\right)\cap\left(\biguplus_{i=1}^{N}\mathcal{S}_{i}\right).
\]
Therefore, we have
\begin{align}
\biguplus_{i=1}^{N}\left(\beta_{i}\mathcal{X}_{o}+t_{i}\right)\cap\mathcal{U}^{i}\subset&\biguplus_{i=1}^{N}\left(\beta_{i}\mathcal{X}_{o}+t_{i}\right)\cap\left(\biguplus_{i=1}^{N}\mathcal{U}^{i}\right)\notag\\
=&(\beta\mathcal{X}_{o}+t)\cap\mathcal{U}_{o}=\mathbb{B}_{n}.\label{eq:InsectMinAssociate}
\end{align}
From the MOA solution, we know $\forall \ i =1,\cdots, N$, $\mathcal{X}^{i}\subset\beta_{i}\mathcal{X}_{o}+t_{i}$.
It follows that
\[
\mathcal{P}^{i}=\biguplus_{i=1}^{N}\mathcal{X}^{i}\cap\mathcal{U}^{i}\subset\biguplus_{i=1}^{N}\left(\beta_{i}\mathcal{X}_{o}+t_{i}\right)\cap\mathcal{U}^{i}.
\]
Combining with~(\ref{eq:InsectMinAssociate}), we have $\mathcal{P}^{i}\subset\mathbb{B}_{n}.$
Hence, $\mathbb{B}_{n}$ is a necessary battery. The battery parameters
$C,\underline{D},\bar{D}$ in~(\ref{eq:NeceBounds2nd}) can be obtained
from formulating the constraint sets on $U$ from $(U-t)/\beta\in\mathcal{X}_{o}$,
while $\underline{E},\bar{E}$ can be obtained by noticing that the Minkowski
sum of hyper-rectangulars can be simply calculated by adding the individual
bounds.\hfill \qed

\bibliographystyle{IEEEtran}

\begin{thebibliography}{10}
\providecommand{\url}[1]{#1}
\csname url@samestyle\endcsname
\providecommand{\newblock}{\relax}
\providecommand{\bibinfo}[2]{#2}
\providecommand{\BIBentrySTDinterwordspacing}{\spaceskip=0pt\relax}
\providecommand{\BIBentryALTinterwordstretchfactor}{4}
\providecommand{\BIBentryALTinterwordspacing}{\spaceskip=\fontdimen2\font plus
\BIBentryALTinterwordstretchfactor\fontdimen3\font minus
  \fontdimen4\font\relax}
\providecommand{\BIBforeignlanguage}[2]{{%
\expandafter\ifx\csname l@#1\endcsname\relax
\typeout{** WARNING: IEEEtran.bst: No hyphenation pattern has been}%
\typeout{** loaded for the language `#1'. Using the pattern for}%
\typeout{** the default language instead.}%
\else
\language=\csname l@#1\endcsname
\fi
#2}}
\providecommand{\BIBdecl}{\relax}
\BIBdecl

\bibitem{Hao2015}
H.~Hao, B.~Sanandaji, K.~Poolla, and T.~Vincent, ``Aggregate flexibility of
  thermostatically controlled loads,'' \emph{{IEEE} Trans. Power Syst.},
  vol.~30, no.~1, pp. 189--198, Jan. 2015.

\bibitem{Hao2015a}
\BIBentryALTinterwordspacing
H.~Hao, A.~Somani, J.~Lian, and T.~E. Carroll, ``Generalized aggregation and
  coordination of residential loads in a smart community,'' in \emph{IEEE Int.
  Conf. Smart Grid Commun.}, May 2015, pp. 67--72. [Online]. Available:
  \url{(Extended Version)
  https://drive.google.com/open?id=0B41ZsRCqdIfYTnlqREE0OFJ4SkU}
\BIBentrySTDinterwordspacing

\bibitem{makarov2009operational}
Y.~Makarov, C.~Loutan, J.~Ma, and P.~de~Mello, ``Operational impacts of wind
  generation on {California} power systems,'' \emph{{IEEE} Trans. Power Syst.},
  vol.~24, no.~2, pp. 1039--1050, May 2009.

\bibitem{helman2010resource}
U.~Helman, ``Resource and transmission planning to achieve a 33\% {RPS} in
  {California}--{ISO} modeling tools and planning framework,'' in \emph{FERC
  Technical Conference on Planning Models and Software}, 2010.

\bibitem{CallawayPIEEE}
D.~S. Callaway and I.~A. Hiskens, ``Achieving controllability of electric
  loads,'' \emph{Proc. {IEEE}}, vol.~99, no.~1, pp. 184 -- 199, 2011.

\bibitem{lu2004state}
N.~Lu and D.~P. Chassin, ``A state-queueing model of thermostatically
  controlled appliances,'' \emph{{IEEE} Trans. Power Syst.}, vol.~19, no.~3,
  pp. 1666--1673, 2004.

\bibitem{Callaway2009}
D.~S. Callaway, ``Tapping the energy storage potential in electric loads to
  deliver load following and regulation with application to wind energy,''
  \emph{Energy Conversion and Management}, vol.~50, no.~5, pp. 1389 --1400, May
  2009.

\bibitem{Mathieu2013}
J.~Mathieu, S.~Koch, and D.~Callaway, ``State estimation and control of
  electric loads to manage real-time energy imbalance,'' \emph{{IEEE} Trans.
  Power Syst.}, vol.~28, no.~1, pp. 430--440, Feb. 2013.

\bibitem{Mathieu2015}
J.~L. Mathieu, M.~Kamgarpour, J.~Lygeros, G.~Andersson, and D.~S. Callaway,
  ``Arbitraging intraday wholesale energy market prices with aggregations of
  thermostatic loads,'' \emph{{IEEE} Trans. Power Syst.}, vol.~30, no.~2, pp.
  763--772, Mar. 2015.

\bibitem{li2015}
S.~Li, W.~Zhang, J.~Lian, and K.~Kalsi, ``Market-based coordination of
  thermostatically controlled loads part {I}: A mechanism design formulation,''
  \emph{{IEEE} Trans. Power Syst.}, vol.~31, no.~2, pp. 1170--1178, 2015.

\bibitem{malhame1985electric}
R.~Malham{\'{e}} and C.-Y. Chong, ``Electric load model synthesis by diffusion
  approximation of a high-order hybrid-state stochastic system,'' \emph{{IEEE}
  Trans. Automat. Contr.}, vol.~30, no.~9, pp. 854--860, 1985.

\bibitem{BashashFathy2013}
S.~Bashash and H.~K. Fathy, ``Modeling and control of aggregate air
  conditioning loads for robust renewable power management,'' \emph{{IEEE}
  Trans. Control Syst. Technol.}, vol.~21, no.~4, pp. 1318--1327, July 2013.

\bibitem{Zhao2015}
L.~Zhao and W.~Zhang, ``A unified stochastic hybrid system approach to
  aggregated load modeling for demand response,'' in \emph{Proc. IEEE Conf.
  Decision and Control (CDC)}, Dec. 2015, pp. 6668--6673.

\bibitem{Zhang2013}
W.~Zhang, J.~Lian, C.-Y. Chang, and K.~Kalsi, ``Aggregated modeling and control
  of air conditioning loads for demand response,'' \emph{{IEEE} Trans. Power
  Syst.}, vol.~28, no.~4, pp. 4655--4664, Nov. 2013.

\bibitem{sanandaji2014fast}
B.~M. Sanandaji, H.~Hao, and K.~Poolla, ``Fast regulation service provision via
  aggregation of thermostatically controlled loads,'' in \emph{Hawaii Int.
  Conf. Syst. Sciences}, 2014, pp. 2388--2397.

\bibitem{Sanandaji2014}
B.~Sanandaji, H.~Hao, K.~Poolla, and T.~Vincent, ``Improved battery models of
  an aggregation of thermostatically controlled loads for frequency
  regulation,'' in \emph{Proc. Amer. Control Conf. (ACC)}, 2014, pp. 38--45.

\bibitem{BarotTaylor2017}
S.~Barot and J.~A. Taylor, ``A concise, approximate representation of a
  collection of loads described by polytopes,'' \emph{International Journal of
  Electrical Power \& Energy Systems}, vol.~84, pp. 55--63, Jan. 2017.

\bibitem{VrettosAndersson2013}
E.~Vrettos and G.~Andersson, ``Combined load frequency control and active
  distribution network management with thermostatically controlled loads,'' in
  \emph{IEEE Smart Grid Commun. 2013 Symp. - Demand Side Management, Demand
  Response, Dynamic Pricing}, 2013.

\bibitem{ZhangShenMathieu2016}
Y.~Zhang, S.~Shen, and J.~Mathieu, ``Distributionally robust chance-constrained
  optimal power flow with uncertain renewables and uncertain reserves provided
  by loads,'' \emph{{IEEE} Trans. Power Syst.}, vol.~32, no.~2, pp. 1378--1388,
  Mar. 2017.

\bibitem{TrovatoTengStrbac2016}
V.~Trovato, F.~Teng, and G.~Strbac, ``Value of thermostatic loads in future
  low-carbon {Great Britain} system,'' in \emph{19th Power Syst. Computation
  Conf. (PSCC)}, June 2016, pp. 1--7.

\bibitem{PJM_M12}
\BIBentryALTinterwordspacing
\emph{M-12: Balancing Operations}, PJM, 2016. [Online]. Available:
  \url{http://www.pjm.com/~/media/documents/manuals/m12.ashx}
\BIBentrySTDinterwordspacing

\bibitem{HughesDominguezGarciaPoolla2016}
J.~T. Hughes, A.~D. Dom{\'{i}}nguez-Garc{\'{i}}a, and K.~Poolla,
  ``Identification of virtual battery models for flexible loads,'' \emph{{IEEE}
  Trans. Power Syst.}, vol.~31, no.~6, pp. 4660--4669, Nov. 2016.

\bibitem{PJMinterconnection}
\BIBentryALTinterwordspacing
{PJM Interconnection}, ``{PJM} regulation data.'' [Online]. Available:
  \url{http://wired.pjm.com/markets-and-operations/ancillary-services.aspx}
\BIBentrySTDinterwordspacing

\bibitem{TrangbaekBendtsen2012}
K.~Trangbaek and J.~Bendtsen, ``Exact constraint aggregation with applications
  to smart grids and resource distribution,'' in \emph{Proc. IEEE Conf.
  Decision and Control (CDC)}, Dec. 2012, pp. 4181--4186.

\bibitem{MuellerSundstroemSzaboEtAl2015}
F.~L. M{\"{u}}ller, O.~Sundstr{\"{o}}m, J.~Szab{\'{o}}, and J.~Lygeros,
  ``Aggregation of energetic flexibility using zonotopes,'' in \emph{Proc. IEEE
  Conf. Decision and Control (CDC)}, Dec. 2015, pp. 6564--6569.

\bibitem{ZhaoHaoZhang2016}
L.~Zhao, H.~Hao, and W.~Zhang, ``Extracting flexibility of heterogeneous
  deferrable loads via polytopic projection approximation,'' in \emph{Proc.
  IEEE Conf. Decision and Control (CDC)}, 2016.

\bibitem{Sanandaji2016}
B.~M. Sanandaji, T.~L. Vincent, and K.~Poolla, ``Ramping rate flexibility of
  residential {HVAC} loads,'' \emph{{IEEE} Trans. Sustain. Energy}, vol.~7,
  no.~2, pp. 865--874, April 2016.

\bibitem{Henk1997}
M.~Henk, J.~Richter-Gebert, and G.~M. Ziegler, ``Basic properties of convex
  polytopes,'' in \emph{Handbook of Discrete and Computational Geometry}, J.~E.
  Goodman and J.~O'Rourke, Eds.\hskip 1em plus 0.5em minus 0.4em\relax Boca
  Raton, FL, USA: CRC Press, Inc., 1997, ch.~15, pp. 243--270.

\bibitem{Schneider1993}
R.~Schneider, \emph{{Convex Bodies: The Brunn-Minkowski Theory}}.\hskip 1em
  plus 0.5em minus 0.4em\relax Cambridge University Press, 1993.

\bibitem{Tiwary2008}
H.~Tiwary, ``On the hardness of computing intersection, union and {Minkowski}
  sum of polytopes,'' \emph{Discrete \& Computational Geometry}, vol.~40,
  no.~3, pp. 469--479, 2008.

\bibitem{Weibel2007}
C.~Weibel, ``{Minkowski} sums of polytopes: Combinatorics and computation,''
  Ph.D. dissertation, {\'{E}}cole Polytechnique F{\'{e}}d{\'{e}}rale De
  Lausanne, 2007.

\bibitem{Eaves1982}
B.~Eaves and R.~Freund, ``Optimal scaling of balls and polyhedra,''
  \emph{Mathematical Programming}, vol.~23, no.~1, pp. 138--147, 1982.

\bibitem{Zhao2016}
L.~Zhao and W.~Zhang, ``A geometric approach to virtual battery modeling of
  thermostatically controlled loads,'' in \emph{Proc. Amer. Control Conf.
  (ACC)}, July 2016, pp. 1452--1457.

\bibitem{Weather}
\BIBentryALTinterwordspacing
``Weather underground: Weather record for {Columbus}.'' [Online]. Available:
  \url{http://www.weatherunderground.com}
\BIBentrySTDinterwordspacing

\bibitem{glpk2006}
\BIBentryALTinterwordspacing
``{GLPK} ({GNU} linear programming kit),'' 2006. [Online]. Available:
  \url{http://www.gnu.org/software/glpk}
\BIBentrySTDinterwordspacing

\bibitem{Lofberg2004}
\BIBentryALTinterwordspacing
J.~L{\"{o}}fberg, ``{YALMIP}: A toolbox for modeling and optimization in
  {MATLAB},'' in \emph{Proc. of the {CACSD} Conf.}, Taipei, Taiwan, 2004.
  [Online]. Available: \url{http://users.isy.liu.se/johanl/yalmip}
\BIBentrySTDinterwordspacing

\bibitem{Mangasarian2002}
O.~Mangasarian, ``Set containment characterization,'' \emph{Journal of Global
  Optimization}, vol.~24, no.~4, pp. 473--480, 2002.

\end{thebibliography}

\begin{IEEEbiography}
    [{\includegraphics[width=1in,height=1.25in,clip,keepaspectratio]{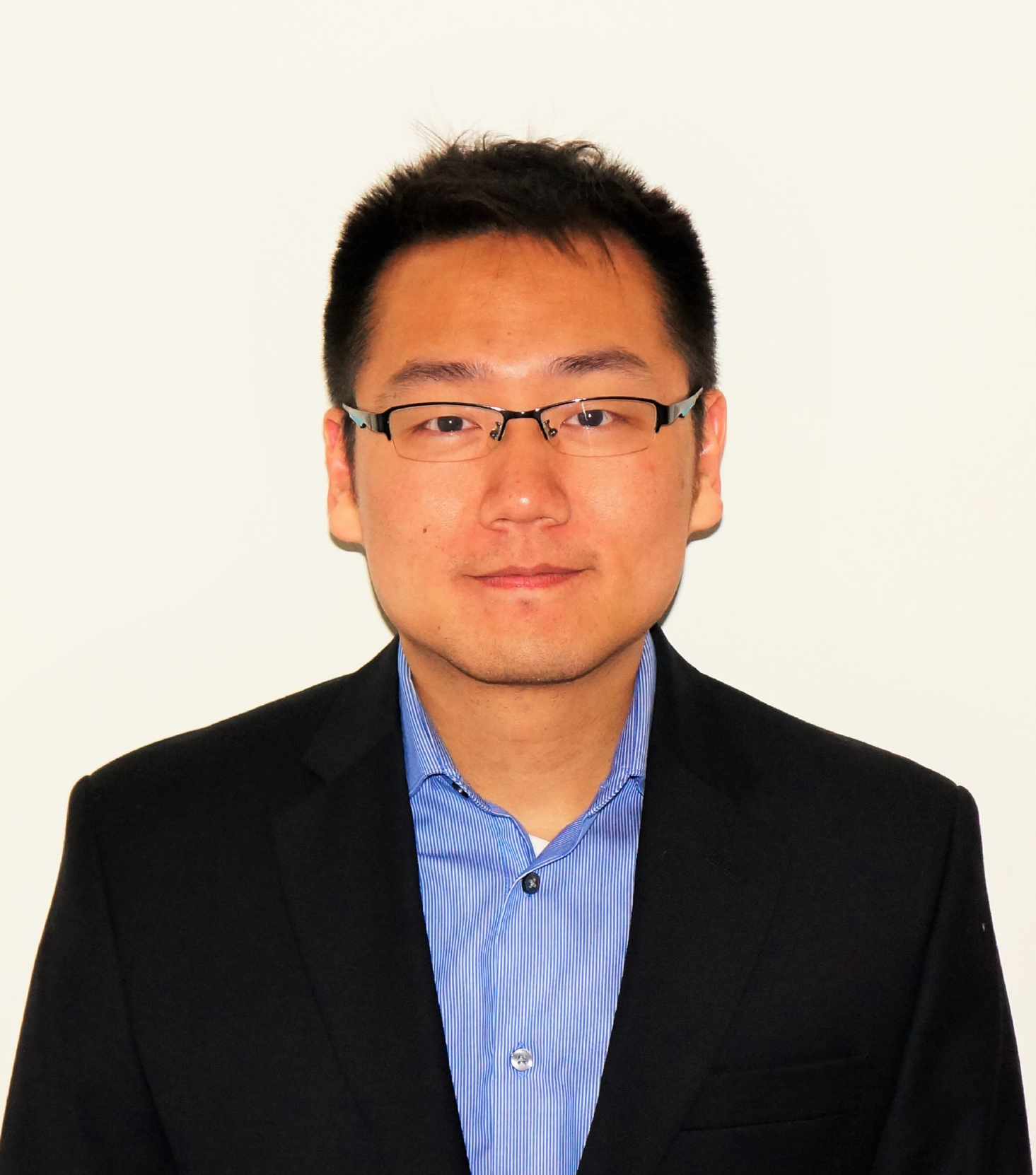}}]{Lin Zhao}
received the B.E. and M.S. degrees in automatic control from the Harbin Institute of Technology, Harbin, China, in 2010 and 2012, respectively. He is currently a Ph.D. student in the Department of Electrical and Computer Engineering, The Ohio State University, Columbus, OH, USA. His current research interests include modeling and control of large-scale complex systems with applications in power systems.
\end{IEEEbiography}

\begin{IEEEbiography}
[{\includegraphics[width=1in,height=1.25in,clip,keepaspectratio]{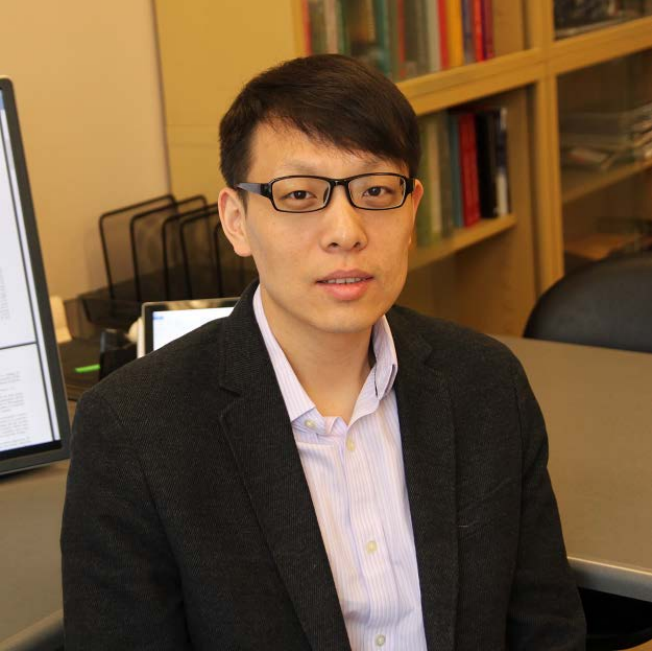}}]{Wei Zhang}
received a B.S. degree in automatic control from the University of Science and Technology of China, Hefei, China, in 2003, and a M.S. degree in statistics and a Ph.D. degree in electrical engineering from Purdue University, West Lafayette, IN, USA, in 2009. Between 2010 and 2011, he was a Postdoctoral Researcher with the Department of Electrical Engineering and Computer Sciences, University of California, Berkeley, CA, USA. He is currently an Assistant Professor in the Department of Electrical and Computer Engineering, Ohio State University, Columbus, OH, USA. His research focuses on control and game theory with applications in power systems, robotics, and intelligent transportations.
\end{IEEEbiography}

\begin{IEEEbiography}
[{\includegraphics[width=1in,height=1.25in,clip,keepaspectratio]{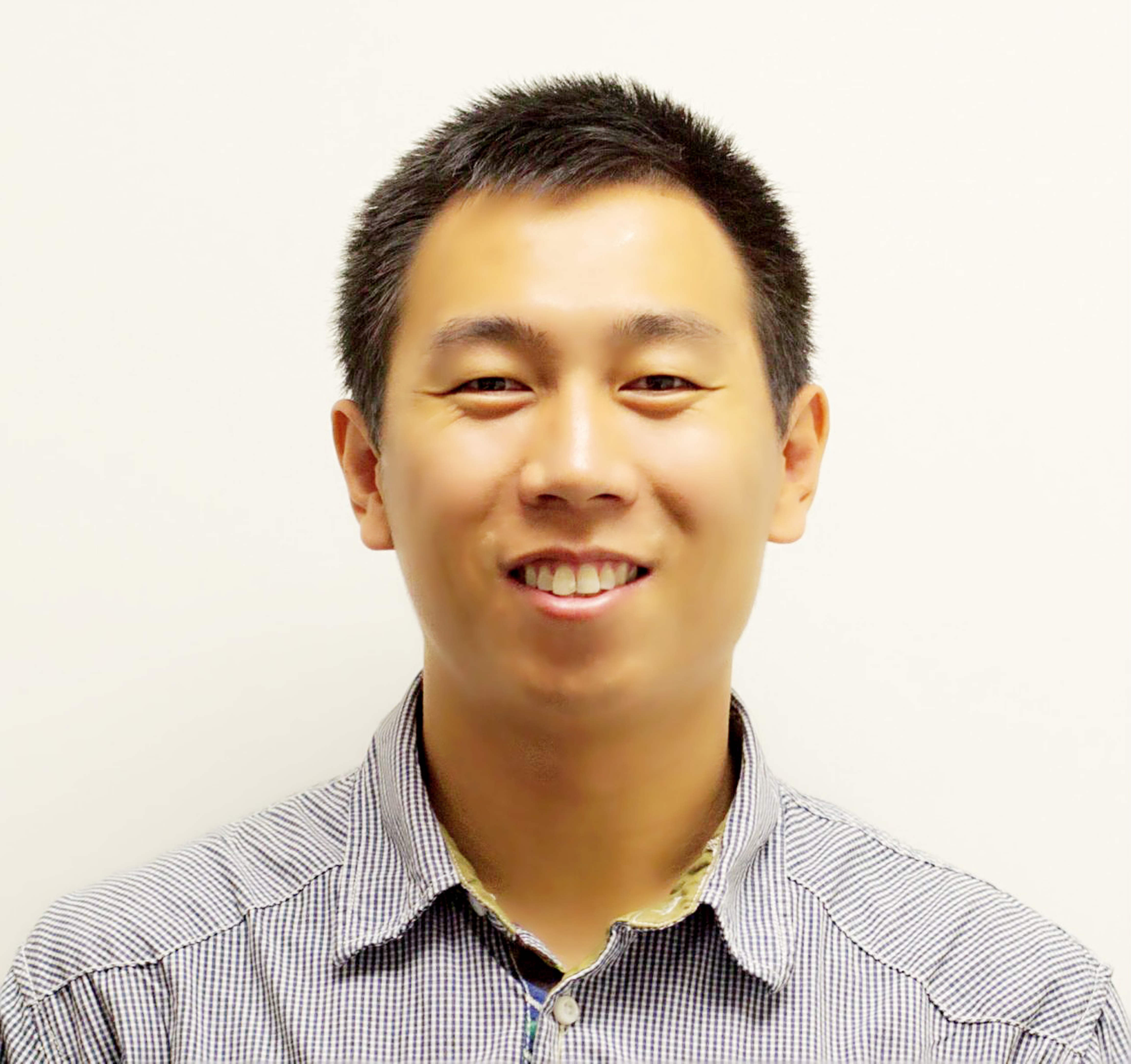}}]{He Hao}
received the B.S. degree in mechanical engineering
and automation from Northeastern University,
Shenyang, China, in 2006, the M.S. degree in
mechanical engineering from Zhejiang University,
Hangzhou, China, in 2008, and the Ph.D. degree
in mechanical engineering from the University of
Florida, Gainesville, FL, USA, in 2012. Between
January 2013 to July 2014, he was a Postdoctoral
Research Fellow with the Department of Electrical
Engineering and Computer Science, University of
California, Berkeley, CA, USA. He is currently a
Staff Scientist in the Electricity Infrastructure and Buildings Division at Pacific
Northwest National Laboratory, Richland, WA, USA. His research interests include
power system, smart buildings, distributed control and optimization, and
large-scale complex systems.
\end{IEEEbiography}

\begin{IEEEbiography}
[{\includegraphics[width=1in,height=1.25in,clip,keepaspectratio]{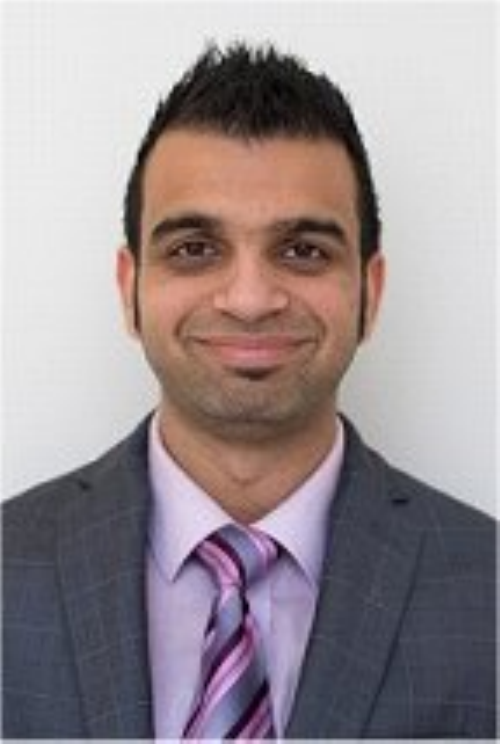}}]{Karanjit Kalsi}
received the M.Eng. degree from the
University of Sheffield, Sheffield, U.K., in 2006, and
the Ph.D. degree in electrical and computer engineering
from Purdue University, West Lafayette, IN,
USA, in 2010. He is currently a Power Systems Research
Engineer at the Pacific Northwest National
Laboratory, Richland, WA, USA.
\end{IEEEbiography}

\end{document}